%% file: Multilevel-Eigensolver.tex
\documentclass[acmtog,copyrightmode]{acmart}

\usepackage{booktabs} 

\usepackage{subfig}
\usepackage{multirow}
\usepackage{wrapfig}
\usepackage{algorithmic}
\usepackage{amsfonts}
\usepackage{amsmath}
\usepackage{soul}

\usepackage[ruled]{algorithm2e} 

\SetAlFnt{\small}
\SetAlCapFnt{\small}
\SetAlCapNameFnt{\small}
\SetAlCapHSkip{0pt}

\acmJournal{TOG}

\begin{document}
\title{The Hierarchical Subspace Iteration Method for Laplace--Beltrami Eigenproblems}

\author{Ahmad Nasikun}
\affiliation{%
  \institution{Delft University of Technology}
  \department{Department of Intelligent Systems}
  \city{Delft}
  \country{The Netherlands}
}
\affiliation{%
  \institution{Universitas Gadjah Mada}
   \department{Department of Electrical and Information Engineering}
  \city{Yogyakarta}
  \country{Indonesia}}
\email{ahmad.nasikun@ugm.ac.id}
\author{Klaus Hildebrandt}
\affiliation{%
  \institution{Delft University of Technology}
  \department{Department of Intelligent Systems}
  \city{Delft}
  \country{The Netherlands}
}
\email{K.A.Hildebrandt@tudelft.nl}
\renewcommand\shortauthors{A. Nasikun \& K. Hildebrandt}

\begin{abstract}
\textit{Sparse eigenproblems are important for various applications in computer graphics. The spectrum and eigenfunctions of the Laplace--Beltrami operator, for example, are fundamental for methods in shape analysis and mesh processing. The Subspace Iteration Method is a robust solver for these problems. In practice, however, Lanczos schemes are often faster. In this paper, we introduce the Hierarchical Subspace Iteration Method (HSIM), a novel solver for sparse eigenproblems that operates on a hierarchy of nested vector spaces. The hierarchy is constructed such that on the coarsest space all eigenpairs can be computed with a dense eigensolver. HSIM uses these eigenpairs as initialization and iterates from coarse to fine over the hierarchy. On each level, subspace iterations, initialized with the solution from the previous level, are used to approximate the eigenpairs. This approach substantially reduces the number of iterations needed on the finest grid compared to the non-hierarchical Subspace Iteration Method. Our experiments show that HSIM can solve Laplace--Beltrami eigenproblems on meshes faster than state-of-the-art methods based on Lanczos iterations, preconditioned conjugate gradients and subspace iterations.}
\end{abstract}

%
%


%
%


\begin{teaserfigure}
\includegraphics[width=\textwidth]{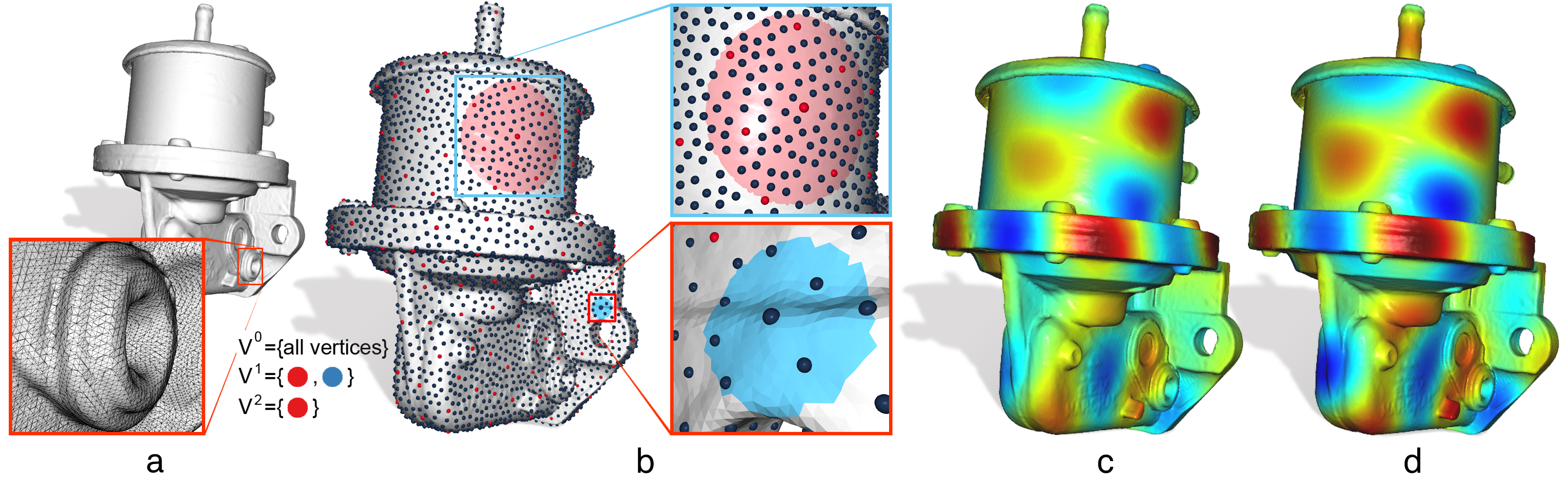}
\caption{
The novel Hierarchical Subspace Iteration Method (HSIM) can efficiently solve eigenvalue problems, such as the computation of the $p$ lowest modes of the discrete Laplace--Beltrami operator on a surface mesh~(a).  
HSIM constructs a set of nested subspaces of the spaces of functions on the mesh by building a vertex hierarchy~(b) and prolongation operators between the levels of the hierarchy. 
HSIM is initialized by solving a dense eigenproblem on the coarsest level, an example of a resulting eigenfunction is shown in~(c).
Then the eigenproblem is solved on each level, from coarse to fine, using subspace iterations initialized with the result from the previous level to finally produce the sought eigenpairs. An example of an eigenfunction is shown in~(d).
}
\label{fig:teaser}
\end{teaserfigure}

\maketitle

\input{body}

\end{document}

%% file: body.tex
\newcommand{\red}[1] {{\color{red}{{#1}}}}
\newcommand{\blue}[1] {{\color{blue}{{#1}}}}
\newcommand{\orange}[1] {{\color{orange}{{#1}}}}
\newcommand{\eg}{\textit{e.g. }}
\newcommand{\ie}{\textit{i.e. }}
\newcommand{\etal}{\textit{et al. }}

\definecolor{revisionColor}{rgb}{0.,0.,0.}
\definecolor{minorRevColor}{rgb}{0.,0.,0.}
\newcommand{\revised}[1]{{\color{revisionColor}{#1}}}
\newcommand{\minorrev}[1]{{\color{minorRevColor}{#1}}}

\newcommand{\level}[4]{{#1}^{#2}_{{#3}_{#4}}}		


\section{Introduction}
Large-scale sparse eigenvalue problems arise in many applications of computer graphics. An important example is the computation of the low and medium frequency spectrum and the corresponding eigenfunctions of the Laplace--Beltrami operator of a surface. These are used in a range of applications in shape analysis and mesh processing. 
Commonly used methods for solving Laplace--Beltrami eigenproblems are based on Lanczos iterations. These are highly efficient solvers for sparse eigenvalue problems. However, in order to be efficient, they combine various extensions of the basic Lanczos iterations, which makes the algorithms complex and introduces parameters that need to be set. 
One problem is that Lanczos iterations are inherently unstable, which can be counteracted by re-starting strategies. Another issue is that Lanczos iterations lead to orthogonal eigenvectors only if the arithmetic is exact. Due to rounding errors, re-orthogonalization strategies are required. 
An alternative to Lanczos schemes is the Subspace Iteration Method (SIM). This method does not suffer from instabilities and is therefore easier to analyze and implement. 
On the other hand, the SIM is often slower than Lanczos schemes. 

In this paper, we introduce the Hierarchical Subspace Iteration Method (HSIM). This method is suitable for computing the eigenpairs in the low and mid frequency part of the spectrum of an operator defined on a mesh, such as the discrete Laplace--Beltrami operator. Our goal is to maintain the benefits of the SIM while reducing the computational cost significantly. One reason why the SIM is expensive is that many iterations are needed before the method converges. Our idea is to take advantage of the fact that low and mid frequency eigenfunctions can be approximated on coarser grids. Instead of working only on the finest grid, we shift iterations to coarser grids. This enables us to perform effective subspace iterations with little computational effort on coarse grids and substantially reduce the number of iterations needed on the finest grid. 

We design the hierarchical solver so that it starts on the coarsest grid. The complexity of this grid is chosen such that the relevant matrices can be represented as dense matrices and all eigenfunctions on the coarsest grid can be efficiently computed with a standard dense eigensolver. Then, the hierarchy is traversed from coarse to fine, whereby the eigenproblem on each grid is solved to the desired accuracy by subspace iterations and the solution on the previous grid is used as an initialization for the subspace iterations. 
To make the subspace iterations more efficient, we use the eigenvalues computed on one grid to determine a value by which we shift the matrix on the next grid. To construct the hierarchy, we use vertex sampling to create a vertex hierarchy and build prolongation operators based on the geodesic vicinity of the samples. 
The prolongation operators are used to define a hierarchy of nested function spaces on the mesh whose degrees of freedom are associated with the vertex hierarchy. 
The advantage of the resulting hierarchy over alternatives, such as mesh coarsening-based hierarchies, is that we obtain a hierarchy of nested spaces. This is of benefit for our purposes because the prolongation to the finer grids then preserves properties of a subspace basis, like its orthonormality.

We evaluate our HSIM scheme on the computation of the lowest $p$ eigenpairs of the Laplace--Beltrami operator, where $p$ ranges from 50 to 5000. Our experiments show that HSIM significantly reduces the number of iterations needed on the finest grid and thus accelerates the SIM method. HSIM has also outperformed three state-of-the-art Lanzcos solvers \revised{and the Locally Optimal Block Preconditioned Conjugate Gradient Method} in our experiments. HSIM was consistently faster than the fastest of the three Lanczos solvers over a range of computations on a variety of meshes and different numbers of eigenpairs to be computed. In particular, for challenging settings, in which more than a thousand eigenpairs needed to be computed, HSIM was up to six times faster than the fastest Lanczos solver. 

\revised{
We expect that applications that need to compute low and medium frequency eigenfunctions of the Laplace-Beltrami operator will benefit from the properties of HSIM, in particular, methods that need to continuously solve new eigenproblems, for example in the context of isospectralization \cite{Cosmo2019,Rampini2019} and geometric deep learning \cite{Bronstein2017}, and methods that need to compute a larger number of eigenfunctions, for example, for shape compression~\cite{Karni2000,Vasa2014}, filtering~\cite{Vallet2008}, and shape signatures~\cite{SOG09}\footnote{
In the supplementary material, we demonstrate that projections into subspaces spanned by Laplace--Beltrami eigenfunctions, and, at the example of the heat kernel signature, that shape signatures can benefit from using a larger number of eigenfunctions.}.
}

\section{Related Work}

\paragraph{Spectral shape analysis and processing}

The eigenfunctions of the Laplace--Beltrami operator on a surface have many properties that make them useful for applications. First, the eigenfunctions form an orthonormal basis of the space of functions on the surface, which generalizes the Fourier basis of planar domains to curved surfaces. With the help of the spectrum and the eigenfunctions, a frequency representation can be associated to functions on a surface and spectral methods from signal and image processing can be generalized to methods for the processing of surfaces. Examples of mesh processing applications that use the Laplace--Beltrami spectrum and eigenfunctions are surface filtering \cite{Vallet2008}, mesh and animation compression \cite{Karni2000,Vasa2014}, quad meshing \cite{Dong06/Brem/Gar/Pas,Huang08Kobbelt,Ling2014}, surface segmentation \cite{Sharma2009,huang09Wicke}, vector field processing \cite{Azencot2013,Brandt2017a}, mesh saliency \cite{Song2014} and shape optimization \cite{Musialski2015}. 
Further properties of the Laplace--Beltrami eigenfunctions are that they are invariant under isometric surface deformation and that they reflect the symmetries of a surface. These properties make them a powerful tool for non-rigid shape analysis. For example, they are used to efficiently compute shape descriptors, such as the the Diffusion 
Distance~\cite{Nadler2005}, the Shape-DNA~\cite{Reuter2005,Reuter2006}, the Global 
Point Signature~\cite{rustamov07}, the Heat Kernel Signature \cite{SOG09}, the Auto Diffusion Function~\cite{Gebal09SGP} and the Wave Kernel Signature
\cite{Aubry2011}. 
\revised{
Moreover the eigenfunctions are the basis for Functional Maps \cite{Ovsjanikov2012,Rustamov2013,Kovnatsky2013,Ovsjanikov2016,Rodola2017,Litany2017}, isospectralization \cite{Cosmo2019,Rampini2019} 
and spectral methods in Geometric Deep Learning \cite{Bruna2014,Boscaini2015,Bronstein2017,Sharp2020}.}

\paragraph{Krylov schemes}

Krylov methods, such as Lanczos schemes for symmetric and Arnoldi schemes for general matrices, are effective solvers for large scale eigenproblems. 
For a comprehensive introduction to Krylov schemes, we refer to the textbook by Saad \shortcite{saad2011numerical}. 
One way to apply Lanczos schemes to generalized eigenproblems, such as the Laplace--Beltrami problem we consider, is to convert them to ordinary eigenproblems by a change of coordinates. 
In particular, if the scalar product is given by a diagonal mass matrix, the change of coordinates is not costly~\cite{Vallet2008}. 
For non-diagonal matrices, the coordinate transformation can be done using a Cholesky decomposition of the mass matrix~\cite{saad2011numerical}.
\textsc{Arpack} \cite{Lehoucq1998} provides implementations of the Implicitly Restarted Lanczos Method for symmetric eigenproblems and the Implicitly Restarted Arnoldi Method for non-symmetric eigenproblems.
\textsc{Arpack} is so widely used that it can be consider to provide reference implementations of the Implicitly Restarted Lanczos and Arnoldi Methods. 
For example, \textsc{Matlab}'s sparse eigensolver \texttt{eigs} interfaces \textsc{Arpack}. 
\textsc{SpectrA} \cite{Qiu2015} is a library offering a C++ implementation of an Implicitly Restarted Lanczos Method build on top of the \textsc{Eigen} matrix library \cite{Guennebaud2010}.   
An alternative to the implicitly restarted Lanczos method is the band-by-band, shift-and-invert Lanczos solver for Laplace--Beltrami eigenproblems on surfaces that was introduced in \cite{Vallet2008}. 

\paragraph{Subspace iterations}
An alternative to Krylov schemes is the Subspace Iteration Method (SIM). It is a robust method for solving generalized sparse eigenproblems and is well-suited for parallelization~\cite{bathe2013subspace}. A comprehensive introduction to the SIM can be found in the textbook by Bathe~\shortcite{bathe2014finite}. 
Matrix shifting is important to make the subspace iterations effective. Different heuristics have been proposed ranging from conservative choices 
\cite{bathe1980accelerated,gong2005comparison} 
to more aggressive shifting strategies \cite{zhao2007accelerated}.   
A recent development is the concept of turning vectors \cite{kim2017bathe} and its extension that includes the turning of turning vectors \cite{wilkins2019e2}. 

\revised{
\paragraph{Preconditioned Eigensolvers}
The lowest eigenpairs of a matrix can be computed by minimizing the Rayleigh coefficient. The Locally Optimal Block Preconditioned Conjugate Gradient Method (LOBPCG) \cite{knyazev2001toward} uses a preconditioned conjugate gradient solver for this minimization. A property of the method is that it does not need to explicitly access the matrix but only needs to evaluate matrix-vector products, which can be of benefit when dealing with large matrices. In recent work \cite{Duersch2018}, an improved basis selection strategy is proposed that improves the robustness of the method when larger numbers of eigenpairs are computed. LOBPCG was used for solving Steklov eigenproblems in~\cite{Wang2019}. 
A method that uses hierarchical preconditioning to approximate a few of the lowest eigenpairs was presented in \cite{Krishnan2013}. 
}
\minorrev{
While LOBPCG is reported to be effective for different eigenproblems, our experiments, see Section~\ref{sect.comparison}, indicate that for the Laplace--Beltrami eigenproblems we consider, HSIM is more efficient. 
}

\paragraph{Approximation schemes}

Schemes for the approximate solution of eigenproblems are static condensation 
\cite{bathe2014finite} in engineering and the Nystr\"om method \cite{Williams2001} 
and random projections \cite{Halko2011} in machine learning.  
Approximation schemes for the Laplace--Beltrami 
eigenproblem on surfaces have been introduced in \cite{Chuang2009,Nasikun2018,Liu2019,Lescoat2020}. 
In contrast to the eigensolvers we consider in this work, these schemes do not provide any guarantee on the approximation quality of the eigenpairs. 

\paragraph{Multigrids on surfaces}

The multigrid hierarchy we need is challenging since we are working with an irregular mesh on a curved surface. 
One way to build a multigrid hierarchy for a surface mesh is to use mesh coarsening algorithms~\cite{Hoppe1996,Aksoylu2005}. 
This is, however, not ideal for our setting because the resulting spaces are not nested, as each space is defined on a different surface. 
Another possibility is to build hierarchical grids on ambient space and then restrict the functions to the surface \cite{Chuang2009}. 
The function spaces generated by this approach, however, do not resemble the linear Lagrange finite elements on the mesh that we want to work with. 
Algebraic multigrids \cite{Stueben2001} are an alternative that would fit our setting. However, unlike the proposed hierarchy, algebraic multigrids only use the operator to build the hierarchy, while we also use the geometry of the surface. 
\revised{A multi-level approach for the computation of the heat kernels on surfaces was introduced in~\cite{Vaxman2010}. In recent work, an intrinsic prolongation operator based on mesh coarsening has been proposed~\cite{Liu2021}.} 

\paragraph{Multilevel eigensolvers}

A traditional multigrid approach to eigenproblems is to treat them as a nonlinear equation and to apply nonlinear multigrid solver to the equation \cite{Hackbusch1979,Brandt1983}. These methods have the advantage that they can be extended or even applied directly to nonlinear eigenproblems. For linear eigenvalue problems, however, this technique is not always efficient because the specific properties of eigenvalue problems are not used when a general nonlinear solver is used. 

Another approach is to integrate a multigrid scheme for solving linear systems into an eigensolver \cite{McCormick1981,Bank1982,Martikainen2001,Arbenz2005}. 
A solver for linear eigenproblems that needs to solve linear systems in every iteration, such as Krylov and subspace iteration methods, is used as an outer iteration. 
In every outer iteration, the linear systems are solved in an inner multigrid loop. 
For our HSIM solver, we use sparse direct solvers for the linear systems, as these are more efficient in our setting than multigrid solvers, see \cite{Botsch2005}. In a different application context, however, it could be useful to use a multigrid linear solver. 

An approach in which also the outer iterations operate on two different grids was proposed in \cite {Xu2001}. In this method, the lowest eigenpair of an elliptic operator is approximated by first computing the eigenpair on the coarse grid and then correcting it by a boundary value problem on the fine grid.  
This two-grid correction scheme was accelerated in \cite{Hu2011} and extended to include matrix shifting in \cite{Yang2011}.  
A multigrid extension of the scheme was introduced in \cite{Lin2015,Chen2016} and later integrated with wavelet bases \cite{Xie2019} and algebraic multigrid procedures \cite{Zhang2015a}. 
The multigrid correction scheme has been used for the computation of Laplace spectra on planar domains~\cite{Hu2011} and parametrized surfaces~\cite{Brannick2015}.
A key difference to HSIM is that HSIM provides users with explicit control of the residual of the resulting eigenpairs. In contrast, the multigrid correction schemes do not provide control over the residual. 
Instead, the resulting residual depends on the approximation quality of the grids in the hierarchy.
We include a discussion and comparison in Section~\ref{sect.comparison}.

\section{Background}

In this section, we first briefly review the Laplace--Beltrami eigenproblem, 
which we use for evaluating the proposed eigensolver. 
Then, we describe the Subspace Iteration Method, which will be the basis of the novel Hierarchical Subspace Iteration Method. 

\subsection{Laplace--Beltrami eigenproblem}

In the continuous case, we consider a compact and smooth surface $\Sigma$ in 
$\mathbb{R}^{3}$. A function $\phi$ is an eigenfunction of the Laplace--Beltrami operator $\Delta$
on $\Sigma$ with eigenvalue $\lambda\in\mathbb{R}$ if%
\begin{equation}
-\Delta\phi=\lambda\phi\label{eq.contEigenProb}
\end{equation}
holds. For discretization, the weak form of (\ref{eq.contEigenProb}) is
helpful. This can be obtained by multiplying both sides of the equation with a
continuously differentiable function $f$ and integrating
\begin{equation}
\int\nolimits_{\Sigma}\text{grad\thinspace}\phi\cdot\text{grad\thinspace
}f\,\text{d}A=\lambda\int\nolimits_{\Sigma}\phi f\,\text{d}%
A. \label{eq.weakEigenProb}%
\end{equation}
On the left-hand side of the equation, we applied integration by parts. 
A function $\phi$ is a solution of (\ref{eq.contEigenProb}) with eigenvalue $\lambda$ if and only if (\ref{eq.weakEigenProb}) holds for all continuously differentiable functions $f$. 
A benefit of the weak form is that evaluating both integrals in
(\ref{eq.weakEigenProb}) only requires functions to be weakly differentiable
(with square-integrable weak derivative) and does not involve differentials of
the surface's metric tensor. 

In the discrete case, $\Sigma$ is a triangle mesh and we consider a
finite-dimensional space of functions defined on the mesh, usually the space
$F$ of continuous functions that are linear polynomials over every triangle.
Then, for functions $\phi,f\in F$, the integrals in (\ref{eq.weakEigenProb})
can be evaluated and $\phi$ is an eigenfunction of the discrete
Laplace--Beltrami operator if there is a $\lambda\in\mathbb{R}$ such that 
(\ref{eq.weakEigenProb}) holds for any $f\in F$. 

Any function in $F$ is uniquely determined by its function values at the
vertices of the mesh. The nodal representation of a function in $F$ is a
vector $\Phi\in \mathbb{R}^{n}$ that lists the function values at all vertices. 
If a nodal vector $\Phi$ is given, the corresponding function in $F$ can be 
constructed by linear interpolation of the function values at the three vertices 
in every triangle. Let $\varphi_{i}\in F$
be the function that takes the value one at vertex $i$ and vanishes at all
other vertices. Then, the stiffness, or cotangent, matrix $S$ and the
mass matrix $M$ are given by
\begin{equation}
S_{ij}=\int\nolimits_{\Sigma}\text{grad\thinspace}\varphi_{i}\cdot
\text{grad\thinspace}\varphi_{j}\text{d}A\qquad\text{and}\qquad M_{ij}
=\int\nolimits_{\Sigma}\varphi_{i}\varphi_{j}\,\text{d}
A.\label{eq.massStiffMatrices}
\end{equation}
Explicit formulas for $S_{ij}$ and $M_{ij}$ can be found in~\cite{Wardetzky2007,Vallet2008}.
The eigenfunctions $\Phi$ and eigenvalue~ $\lambda$ can be computed as the
solution to the eigenvalue problem
\begin{equation}
S\,\Phi=\lambda M\,\Phi.\label{eq.eigenProb}%
\end{equation}
This is a sparse, generalized eigenvalue problem where $M$ is symmetric and
positive definite and $S$ is symmetric.
We refer to \cite{Hildebrandt2006,Crane2013} for more background on the discretization of the Laplace--Beltrami operator on surfaces.  

\subsection{Subspace iteration method}

The subspace iteration method (SIM) is an approach for computing 
eigenpairs of generalized eigenvalue problems such as~(\ref{eq.eigenProb}).
We outline SIM in Algorithm~\ref{alg:basic_SIM}. The input to the method
are the stiffness and mass matrices $S,M\in\mathbb{R}^{n\times n}$, 
a matrix $\Phi\in\mathbb{R}^{n\times q}$ that specifies an initial subspace 
basis, the number of desired eigenpairs $p$, a tolerance $\varepsilon$ and a 
shifting value $\mu$. The dimension $q$ of the subspace needs to be 
larger or equal to $p$.
We will first discuss the subspace iterations without shifting, \textit{i.e.} assuming $\mu=0$, and then 
discuss choices of convergence test, subspace dimension, initial subspace basis, 
shifting value and linear solver. 

SIM iteratively modifies the initial basis, which consists of $q$ vectors $\Phi_i$, such that it more and more
becomes the desired eigenbasis. In each iteration, first an inverse	
iteration is applied to all $q$ vectors (Algorithm \ref{alg:basic_SIM}, line
4), thereby increasing the low-frequency components in the vectors. For this,
$q$ linear systems of the form
\begin{equation}
(S-\mu M)\Psi_{i}=M\Phi_{i}%
\end{equation}
need to be solved. 
The second step in each iteration is to
solve the eigenproblem restricted to the subspace spanned by the vectors
$\Psi_{i}$ (Algorithm \ref{alg:basic_SIM}, lines 5--7). For this, the reduced
stiffness and mass matrices are computed and the $q$-dimensional dense
eigenproblem is solved using a dense eigensolver, \emph{e.g.} based on a QR factorization. The third step is to replace the 
current subspace basis with the eigenbasis (Algorithm \ref{alg:basic_SIM}, line 8). 
The inverse iterations amplify the low frequencies in the subspace 
basis. The second and third steps are needed in order to prevent the vectors 
from becoming linearly dependent. Without these steps, the vectors would all converge to the lowest eigenvector. 

\begin{algorithm}[tbh]
\DontPrintSemicolon
\LinesNumbered
\KwIn{Stiffness matrix $S\in\mathbb{R}^{n\times n}$, mass matrix $M\in\mathbb{R}^{n\times n}$, initial vectors $\Phi\in\mathbb{R}^{n\times q}$, number of eigenpairs $p$, tolerance $\varepsilon$, shifting value $\mu$}
\KwOut{Matrix $\bar\Lambda$ with lowest eigenvalues of (\ref{eq.eigenProb}) on diagonal and $\Phi$ listing eigenvectors as columns. First $p$ pairs converged.}
\SetKwFunction{FMain}{SIM}
\SetKwProg{Fn}{Function}{:}{}
\Fn{\FMain{$S, M, \Phi, p, \varepsilon, \mu$}}{
	Compute sparse factorization:  $LDL^T=S-\mu M$\;
	\Repeat{pairs ($\bar\Lambda_{ii},\Phi_i$) pass convergence test (\ref{eq.conv2}) for all $i\leq p$}{
		Solve using factorization: $(S-\mu M)\Psi = M \Phi$\;
		Compute reduced stiffness matrix: $\bar{S}\gets\Psi^{T} S \Psi$\;
		Compute reduced mass matrix: $\bar{M}\gets\Psi^{T} M \Psi$\;
		Solve dense eigenproblem: $\bar{S} \bar{\Phi} =  \bar{M} \bar{\Phi} {\bar\Lambda}$\;
		Update vectors: $\Phi\gets \Psi\bar{\Phi}$\; 
	}
	\Return $\bar\Lambda$ and $\Phi$
}
\textbf{End Function}
\caption{Subspace Iteration Method}
\label{alg:basic_SIM}
\end{algorithm}

\paragraph{Convergence test} 
The final step of each iteration is the convergence check, which tests whether
or not the first $p$ eigenvectors have converged. For each eigenpair $\Phi
_{i}$ and $\lambda_{i}$, the \minorrev{relative $M^{-1}$-norm\footnote{\minorrev{The $M^{-1}$-norm
is given by ${\left\Vert S\Phi\right\Vert _{M^{-1}}=}%
\sqrt{(S\Phi)^{T}M^{-1}S\Phi}.$ The reason, we use the $M^{-1}$-norm is that
$S\Phi$ is an integrated quantity and $M^{-1}S\Phi$ is the corresponding
function (pointwise quantity). The $M$-norm, which is the discrete $L^{2}%
$-norm, of the pointwise quantity is the same as the $M^{-1}$-norm of the
integrated quatity, ${\left\Vert S\Phi\right\Vert _{M^{-1}}=}\sqrt{(S\Phi
)^{T}M^{-1}S\Phi}=\sqrt{(M^{-1}S\Phi)^{T}MM^{-1}S\Phi}{=\left\Vert M^{-1}%
S\Phi\right\Vert _{M}}$. For more background, we refer to \cite{Wardetzky2007}.}} 
of the residual of equation
(\ref{eq.eigenProb}) is computed
\begin{equation}
{\frac{{\left\Vert {S\Phi_{i}-\lambda_{i}M\Phi_{i}}\right\Vert_{M^{-1}} }%
}{{\left\Vert S\Phi_{i}\right\Vert_{M^{-1}} }}}<\varepsilon\label{eq.conv2}%
\end{equation}
and the test is passed if it is below the threshold $\varepsilon$.}
The choice of the value for the convergence tolerance depends on the application 
context. 
In most of our experiments, we used $\varepsilon=10^{-2}$, \revised{which based on our experiments, see Section \ref{sect.experiments}, we consider appropriate for applications in shape analysis and spectral mesh processing. } 

\paragraph{Subspace dimension} 
The choice of the dimension $q$ affects the computational cost per iteration 
and the number of iterations needed for convergence. 
A larger subspace size increases the computational cost per iteration as 
more linear systems have to be solved (line 4 of Algorithm \ref{alg:basic_SIM}) 
and the dimension of the dense eigenproblem (line 7) increases. On the other 
hand, the algorithm terminates when the subspace contains (good enough 
approximations of) the lowest $p$ eigenvectors. This is easier to achieve if the 
subspace is larger. Therefore, with a larger subspace, fewer iterations may be 
needed. 
It is suggested to set $q=\max\{2p,p+8\}$ in \cite{
bathe2013subspace}. In our experiments, we found $q=\max\{1.5p,p+8\}$
to be more efficient for the eigenproblems we consider. 

\paragraph{Initialization}
The subspace basis $\Phi$ can be initialized with a random matrix. 
An alternative is to use information extracted from the matrices for initialization, which can help to reduce the required number of subspace iterations. 
One heuristic from \cite{bathe2014finite} is to use the diagonal of
the mass matrix $M$ as 
the first column of the matrix representing initial vectors, 
random entries for the last column, and unit vectors $e_i$ 
with entry $+1$ at the degree of freedom with the smallest ratio of $k_{ii}/m_{ii}$
for the remaining $q-2$ columns.

\paragraph{Shifting}
One way to make the subspace iterations more effective is to shift the matrix 
$S$, which means to replace it with the shifted matrix $S-\mu M$. The shifted 
matrix keeps the same eigenvectors while the eigenvalues are shifted by $-\mu$. 
As a consequence, the inverse iteration, line 4 of Algorithm 
\ref{alg:basic_SIM}, focuses on enhancing the frequencies around $\mu$ instead of 
around zero. This can help to reduce the number of iterations required for convergence. 
Different heuristics for setting the shifting value have been proposed. A conservative choice is to set $\mu$ to the average of the last two converged eigenvalues~\cite{bathe1980accelerated}. Alternative shifting strategies are to set $\mu$ to the average of the last converged and the first non-converged eigenvalue~\cite{wilson1983eigensolution} or to the average of the first two non-converged eigenvalues~\cite{gong2005comparison}. An aggressive shifting technique that places $\mu$ further into the range of the non-converged eigenvalues is shown to accelerate the SIM in \cite{zhao2007accelerated}.

\paragraph{Direct solver}
For the inverse iterations of the subspace basis, line~4 of 
Algorithm~\ref{alg:basic_SIM}, $q$ linear systems with the same matrix $S-\mu M$
need to be solved. 
It can be effective to use a direct solver for this task since a factorization 
once computed can be used to solve all the systems. 
Since the shifted matrix is not positive definite, we use a sparse symmetric indefinite decomposition $LDL^{T}=S-\mu M$. 


\section{Hierarchical Subspace Iteration Method}

In this section, we introduce the Hierarchical Subspace Iteration Method
(HSIM). We first describe the construction of the hierarchy of function spaces
on a mesh. Then, we detail the multilevel eigensolver that operates on the hierarchy.

\subsection{Hierarchy construction}

\label{ssec.hierarchy}Important goals for the construction of the hierarchy
are that the construction is fast since the hierarchy must be built as part
of the HSIM algorithm, that the basis functions are locally supported and the
prolongation and restriction operators are sparse, and that the functions
spaces are nested. Moreover, the function spaces need to be able to approximate 
low and mid frequency functions well. 

We describe the construction of the subspaces in three steps. First, we
describe the construction of a hierarchy on the set of vertices of the mesh.
Then, we define prolongation and restriction operators that act between the
levels of the vertex hierarchy. Finally, we explain how the vertex hierarchy
and the operators can be used to obtain the hierarchy of nested function spaces.

\begin{algorithm}[tbh]
\DontPrintSemicolon
\LinesNumbered
\KwIn{Surface mesh $\Sigma$, number of levels $T$, number of vertices per level $n^1,n^2,\dots,n^{T-1}$}
\KwOut{Sets of vertex indices $V^1,V^2,\dots,V^{T-1}$}
	$V^T \gets \{$Random number from $\{0,1,\dots,|V_\Sigma|-1 \}\}$\;
	$\tau\gets T-1$\;
	\Repeat{$\tau=0$}{
		$V^{\tau}\gets V^{\tau+1}$\;
		\Repeat{$|V^\tau|=n^\tau$}{
			$V^{\tau}\gets V^{\tau}\cup\{$Index of vertex farthest away from $V^{\tau}\}$\;
		}
		$\tau \gets \tau-1$\;
	}
	\Return $V^1,V^2,\dots,V^{T-1}$
\caption{Construction of the vertex hierarchy}
\label{alg:vertexHierarchy}
\end{algorithm}

\paragraph{Vertex hierarchy}

We consider a hierarchy with $T$ levels ranging from $0$ to $T-1$, where $0$
is the finest level. We denote by $V^{\tau}$ the set of vertices in level
$\tau$ and by $n^{\tau}$ the number of vertices in $V^{\tau}$. The sets
$V^{\tau}$ are nested, $\level{V}{\tau}{}{}\subset\level{V}{\tau-1}{}{}$, and
$\level{V}{0}{}{}$ is the set of all vertices of the mesh. Since we will solve
a dense eigenproblem to get all eigenpairs at the coarsest level, we want to
control the number $n^{T-1}$ of vertices in $\level{V}{T-1}{}{}$, which we set
to
\begin{equation}\label{eq.nT-1}
n^{T-1}=\text{max}\{\left\lceil 1.5 p \right\rceil,1000\}.
\end{equation}
The numbers of vertices in the other levels are determined by the growth rate
$\mu$%
\begin{equation}
n^{\tau}=\mu\,n^{\tau+1},
\end{equation}
where $\mu$ is given by 
\begin{equation}\label{eq.growthRate}
\mu=\sqrt[T]{\frac{n^{0}}{n^{T-1}}}.
\end{equation}
The trade-off for the choice of the number of levels is that a larger number of levels helps to reduce the required number of iterations on the finest level. On the other hand,  each level adds computational cost, \textit{e.g.} for computing the reduced matrices $S^{\tau}$ and $M^{\tau}$. In our experiments, we found HSIM to be most effective with a low number of levels. We used three levels in most cases and opted for two levels when only a small number of eigenpairs, \ie $p\leq200$, needs to be computed.

To form the sets $V^{\tau}$, we use a scheme based on farthest point sampling
\cite{Eldar1997}. The set $V^{T-1}$ is initialized to contain one random
vertex. Then, iteratively the vertex farthest away from all the vertices that
are already in $V^{T-1}$ is added to $V^{T-1}$ until the desired number of
vertices is reached. The sets $V^{T-2}$ to $V^{1}$ are created in a similar
manner. The scheme is summarized in Algorithm~\ref{alg:vertexHierarchy}. 
Most expensive in this algorithm is the repeated computation of the farthest points  (line~6). 
These computations can be accelerated by maintaining a distance field
that stores for each vertex of the mesh the distance to the closest vertex in the
current set $V^{\tau}$. Since the vertices are inserted one after another, in
each iteration the distance field only needs to be updated locally around the newly
inserted vertex, and the maximum of the field has to be computed. We compute
the distances between vertices using Dijkstra's algorithm on the edge graph
with weights corresponding to the length of the edges. We found Dijkstra's
distance a sufficient approximation of the geodesic distance for our purposes in our
experiments. Alternatively, the Short-Term Vector Dijkstra (STVD) algorithm~\cite{Campen2013} 
could be used, which computes a more accurate approximation of the geodesic
distance while still keeping computations localized.
\revised{The supplementary material includes examples that illustrate that farthest point sampling generates hierarchies that are suitable for our purposes.}

\begin{figure}[t]
  \centering
  \includegraphics[width=1.0\linewidth]{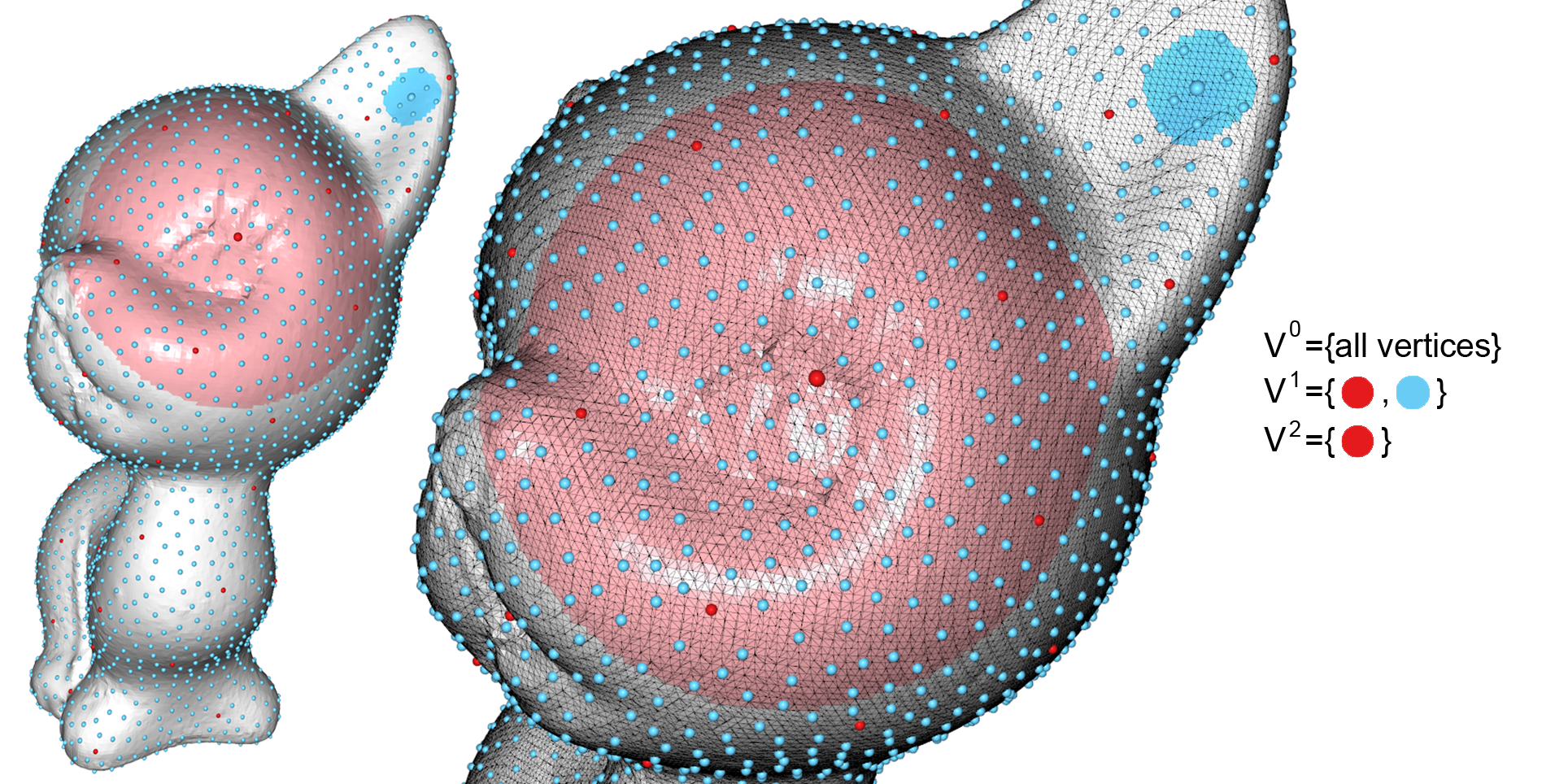}
  \parbox[t]{1.0\columnwidth}{\relax  }
  \caption{\label{fig:hierarchy}
	Illustration of a vertex hierarchy with three sets $V^2\subset V^1 \subset V^0$. The coarsest set $V^2$ consists of the red vertices, $V^1$ of the blue and the red vertices and $V^0$ of all vertices of the mesh. The light red and light blue regions are geodesic disks of radii $\rho^{2}$ and $\rho^1$ around the highlighted red and blue vertices in the centers of the regions and illustrate the support regions of the highlighted vertices.
}
\end{figure}

\paragraph*{Prolongation and restriction}

A function on level $\tau$ is represented by a vector $f^{\tau}%
\in\mathbb{R}^{n^{\tau}}$. We will first describe the prolongation and
restriction operators and show in the next paragraph how the prolongation
operator can be used to construct the piecewise linear polynomial
corresponding to a vector $f^{\tau}$. The $\tau^{th}$ prolongation operator is
given by a matrix $U^{\tau}\in\mathbb{R}^{\level{n}{\tau}{}{}\times
{\level{n}{\tau+1}{}{}}}$ that maps vectors $f^{\tau+1}\in\mathbb{R}%
^{n^{\tau+1}}$ representing functions on level $\tau+1$ to vectors $f^{\tau
}\in\mathbb{R}^{n^{\tau}}$ representing functions on the finer level $\tau$.
The restriction operator maps from level $\tau$ to the coarser level $\tau+1$
and is given by the transpose $U^{\tau T}$ of the prolongation matrix. This
relationship of the prolongation and restriction operators ensures that the
restricted matrices $S^{\tau}$ and $M^{\tau}$, see lines 4 and 5 of
Algorithm~\ref{alg:ours} for a definition of the matrices, on all levels are symmetric.

The $i^{th}$ row of $U^{\tau}$ describes how the value associated with the
$i^{th}$ vertex of level $\tau+1$ is distributed among the vertices on level
$\tau$. This means that the entry $U_{ij}^{\tau}$ is a weight describing how strongly the vertex $j$ on level $\tau$ is influenced by the vertex $i$ on level $\tau+1$ during the prolongation. This weight decreases with increasing geodesic distance of the
vertices. To obtain sparse operators, the weight vanishes when the distance of
the vertices reaches a threshold $\rho^{\tau}$, which differs per level. We
set $\rho^{\tau}$ to be
\begin{equation}
\level{\rho}{\tau}{}{}=\sqrt{\frac{\sigma A}{\level{n}{\tau}{}{}\pi}%
}.\label{eq.radius}%
\end{equation}
where $A$ is the area of the surface and $\sigma$ is a control parameter. This
choice of $\rho^{\tau}$ yields matrices $U^{\tau}$ that have about $\sigma$
non-zero entries per row. For our experiments, we choose $\sigma=7$. The
reasoning behind (\ref{eq.radius}) is that we want the sum of the areas of the
geodesics disks of radius $\rho^{\tau}$ around all the vertices of level
$\tau$ to be $\sigma$ times the area of the surface. To make this idea easily
computable, we replace the combined areas of all the geodesic disks by $n^{\tau
}$ times the area of the Euclidean disk of radius $\rho^{\tau}$. 

To construct the matrices $U^{\tau}$, we first construct preliminary matrices
${\level{\tilde{U}}{\tau}{}{}}\in\mathbb{R}^{\level{n}{\tau}{}{}\times
{\level{n}{\tau+1}{}{}}}$ that have the entries%
\begin{equation}\label{eq.basis}
\tilde{U}_{ij}^{\tau}=\left\{
\begin{tabular}
[c]{cc}%
$1-\frac{d(v_{i}^{\tau+1},v_{j}^{\tau})}{\rho^{\tau}}$ & for~$d(v_{i}^{\tau+1},v_{j}^{\tau}%
)\leq\rho^{\tau}$ \\
$0$ & for $d(v_{i}^{\tau+1},v_{j}^{\tau})>\rho^{\tau}$%
\end{tabular}
\ ,\right.
\end{equation}
where $d(v_{i}^{\tau+1},v_{j}^{\tau})$ is the geodesic distance of the
$i^{th}$ vertex of $V^{\tau+1}$ to the $j^{th}$ vertex of $V^{\tau}$. 
The matrix $\level{U}{\tau}{}{}$ is then obtained by normalizing the rows of
${\level{\tilde{U}}{\tau}{}{}}$
\begin{equation}
U_{ij}^{\tau}=\frac{1}{\sum_{j=1}^{n^{\tau}}\tilde{U}_{ij}^{\tau}}\tilde
{U}_{ij}^{\tau}.
\end{equation}
The normalization ensures that all function spaces will include the constant
functions. This is of benefit for our purposes as the constant functions make
up the kernel of the Laplace--Beltrami operator. Another property is that the
set of functions on each level forms a partition of unity. 
\revised{As for the sampling scheme, we use Dijkstra's distance on the weighted edge graph of the mesh in our experiments to approximate the geodesic distance. A discussion of two alternatives, the Short-Term Vector Dijkstra algorithm~\cite{Campen2013} and the Heat Method~\cite{crane2013geodesics}, is included to the supplementary material.  
}

\paragraph*{Function spaces}

So far we have considered abstract vectors $f^{\tau}\in\mathbb{R}^{n^{\tau}}$. 
Now, we describe how the continuous piecewise linear polynomial
corresponding to $f^{\tau}$ can be constructed. On the finest level, any
$f^{0}\in\mathbb{R}^{n^{0}}$ is the nodal vector, which lists the function values of the continuous, piecewise linear polynomial at the vertices. To get the
continuous, piecewise linear polynomial that corresponds to a 
$f^{\tau}\in\mathbb{R}^{n^{\tau}}$ for any $\tau$, we use the prolongation operators to lift $f^{\tau}$ to the finest level. The resulting vector%
\begin{equation}
U^{0}U^{1}...U^{\tau-1}f^{\tau}\label{eq.liftingFunctions}%
\end{equation}
is the nodal vector of the continuous, piecewise linear polynomial
corresponding to $f^{\tau}$. 
By construction, the resulting function spaces are nested and the functions are locally supported. 
The HSIM algorithm does not need to lift the functions using
(\ref{eq.liftingFunctions}). Instead, the reduced stiffness and mass matrices
$S^\tau$ and $M^\tau$ are directly computed for each level.

\subsection{Hierarchical Solver}
The HSIM is outlined in Algorithm~\ref{alg:ours}. 
The algorithm starts with preparing the multilevel subspace iterations. 
First, the number of levels, the vertex hierarchy and the prolongation matrices $U^{\tau}$ are computed. Then the reduced stiffness and mass matrices,
$S^{\tau}$ and $M^{\tau},$ for all levels are constructed from fine to coarse
starting with level 1. In this computation, we benefit from the fact that the 
prolongation matrices $U^{\tau}$ are highly sparse. The next step is to determine the dimension $q$ of the subspace that is used. 
Our experiments indicate that values between $q=1.5p$ and $q=2p$ are suitable. 
Following \cite{bathe2013subspace}, we set $q=p+8$ for small values of $p$
\begin{equation}
q=\max\{\left\lceil 1.5p\right\rceil ,p+8\}.
\end{equation} 
The last step before the multilevel iterations start is the computation of 
an initial subspace. 
This is done by solving the eigenproblem on the coarsest level of the hierarchy 
completely using a dense eigensolver. 
The dimension of the coarsest space is  
chosen, see (\ref{eq.nT-1}), such that the dense eigenproblem can be solved 
efficiently.

\begin{algorithm}[tbh]
\DontPrintSemicolon
\LinesNumbered
\KwIn{Stiffness and mass matrices of finest level $S^0,M^0\in\mathbb{R}^{n\times n}$, number of eigenpairs $p$, number of levels $T$, tolerance $\varepsilon$}
\KwOut{$p$ lowest eigenpairs of the generalized eigenproblem (\ref{eq.eigenProb})}
\SetKwFunction{FMain}{HSIM}
\SetKwProg{Fn}{Function}{:}{}
\Fn{\FMain{$S^0, M^0, p, T, \varepsilon$}}{
	Compute vertex hierarchy (Section~\ref{ssec.hierarchy}) \;
	Build matrices $\level{U}{\tau}{}{}$ for $\tau=0,1,\dots,T-2$ (Section \ref{ssec.hierarchy}) \;
	\For{$\tau \gets 1$ to $T-1$}{
		Build level $\tau$ stiffness matrix: $S^{\tau}\gets(U^{\tau -1})^T S^{\tau -1} U^{\tau -1}$\;
		Build level $\tau$ mass matrix: $M^{\tau}\gets(U^{\tau -1})^T M^{\tau -1} U^{\tau -1}$\;
	}
	Set size of subspace: $q \gets$ max$(\left\lceil 1.5p\right\rceil,p+8)$\;
	Compute first $q$ eigenpairs of $S^{T-1} \Phi^{T-1} = \Lambda^{T-1} M^{T-1} \Phi^{T-1}$ \;
	\For{$\tau \gets (T-2)$ to $0$}{
		Prolongation of subspace basis: $\Phi^\tau \gets {\level{U}{\tau}{}{}} \Phi^{\tau+1} $\;
		Set shifting parameter: $\mu \gets \Lambda^{\tau+1}_{jj}$ with $j={\lfloor \frac{p}{10} \rfloor}$\;
		$(\Lambda^\tau,\Phi^\tau)\gets$ \textbf{SIM}($S^\tau, M^\tau, \Phi^\tau, p,\varepsilon, \mu$) \;
	}
	\Return First $p$ diagonal entries of $\Lambda^0$ and first $p$ columns of $\Phi^0$
}
\textbf{End Function}
\caption{Hierarchical Subspace Iteration Method}
\label{alg:ours}
\end{algorithm}

The multilevel
iterations traverse the hierarchy from coarse to fine starting with the second
coarsest level. At each level, the eigenproblem is solved up to the 
tolerance by subspace iterations. The subspace iterations are initialized with
the eigenvectors computed at the coarser level. To make the subspace iteration
more effective, we use the approximate eigenvalues computed on the previous
level to specify a shifting parameter for the iterations on the current level. We 
employ an aggressive shifting strategy, 
which sets the shifting
value to be the estimated eigenvalue with index $\left\lfloor p/10\right\rfloor
$. The shifting value is set only once for each level and used for all subspace 
iterations on this level. Then, only one sparse
factorization of the shifted stiffness matrix $S^{\tau}-\mu^{\tau}M^{\tau}$
has to be computed per level. This way we achieve that, on the one hand, the
shifting value is regularly updated, while, on the other hand, no additional
factorizations have to be computed.

To further accelerate the subspace iterations, we do not perform additional inverse iterations, step 4 of the Algorithm \ref{alg:basic_SIM}, on the lowest $r$ vectors that are already converged. However, to avoid error accumulation, we stop the iteration of vectors only after their residual, eq. (\ref{eq.conv2}), has reached one tenth of the specified tolerance $\varepsilon$.  
A further acceleration is achieved by performing two inverse iterations before orthonormalizing the vectors. Thus we execute step~4 of Algorithm~\ref{alg:basic_SIM} twice before we continue with step~5.  

The subspace iteration method converges quickly when the desired 
eigenspace is close to the initial subspace. Our hierarchical method 
makes use of this property by providing the subspace iterations on each 
level with the solution from the coarser level. As a result, only few 
iterations are needed on each level. 
In particular, the multilevel strategy substantially reduces the necessary number of iterations on the finest level compared to SIM.
The price to pay is that the hierarchy has 
to be built and iterations on the coarse levels are needed. 
Nevertheless, HSIM is about 4-8 times faster than SIM in our experiments. The highest acceleration is achieved in the difficult case that a large number of eigenvectors must be computed.

\begin{figure}[b]
  \centering
  \includegraphics[width=1.0\linewidth]{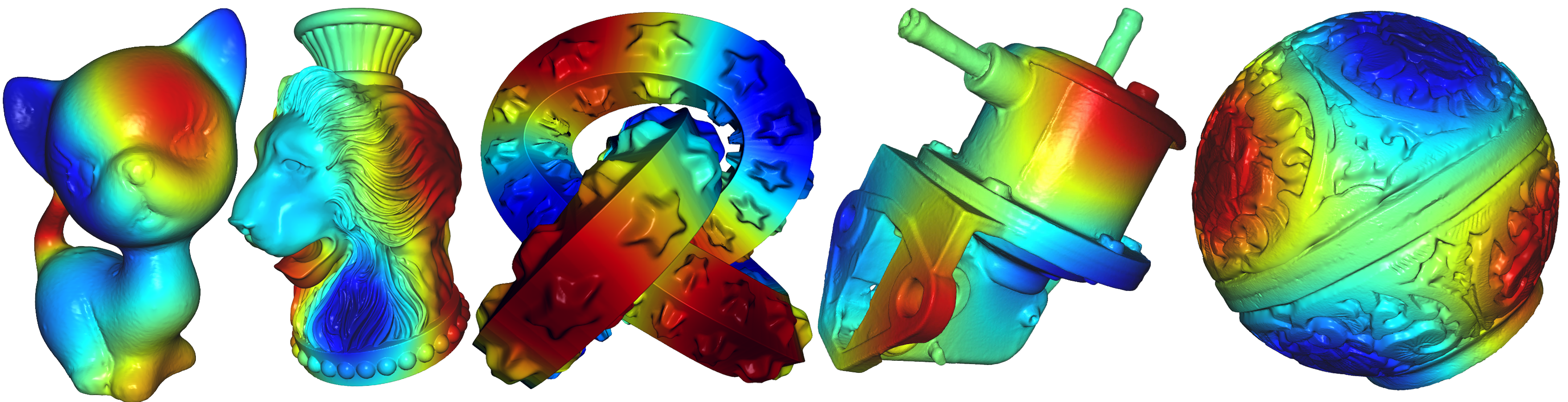}  
  \parbox[t]{1.0\columnwidth}{\relax  }
  \caption{\label{fig:eigenfunction}
  \revised{The 10th eigenfunction of the Laplace-Beltrami operator of the models listed in Table~\ref{tab:construction}.}
}
\end{figure}
\section{Experiments}
\label{sect.experiments}

\paragraph*{Implementation}
Our implementation of HSIM uses \textsc{Eigen} \cite{Guennebaud2010} for linear algebra functionalities
and \textsc{LibIGL} \cite{Jacobson16libigl} for geometry processing tasks. 
\textsc{OpenMP} is used to solve the linear systems in each subspace iteration, step 4 of Algorithm~\ref{alg:basic_SIM}, in parallel and to compute the prolongation and projection matrices in parallel during hierarchy construction. 
Moreover, we solve the low-dimensional eigenproblems at the coarsest level of the hierarchy, step 9 of Algorithm~\ref{alg:ours}, and in each subspace iteration, step 7 of Algorithm~\ref{alg:basic_SIM}, on the GPU using a direct solver for dense generalized eigenproblems from \textsc{CUDA}'s \textsc{cuSolver} library. 
\minorrev{For our experiments, we used an Alienware Area-51 R3 home desktop with a AMD Ryzen Threadripper 1950x (16 core) processor and 24GB of RAM, equipped with NVIDIA GeForce GTX 1080 Ti graphics card with 11GB memory.}

\begin{figure}[b]
  \centering
  \includegraphics[width=1.0\linewidth]{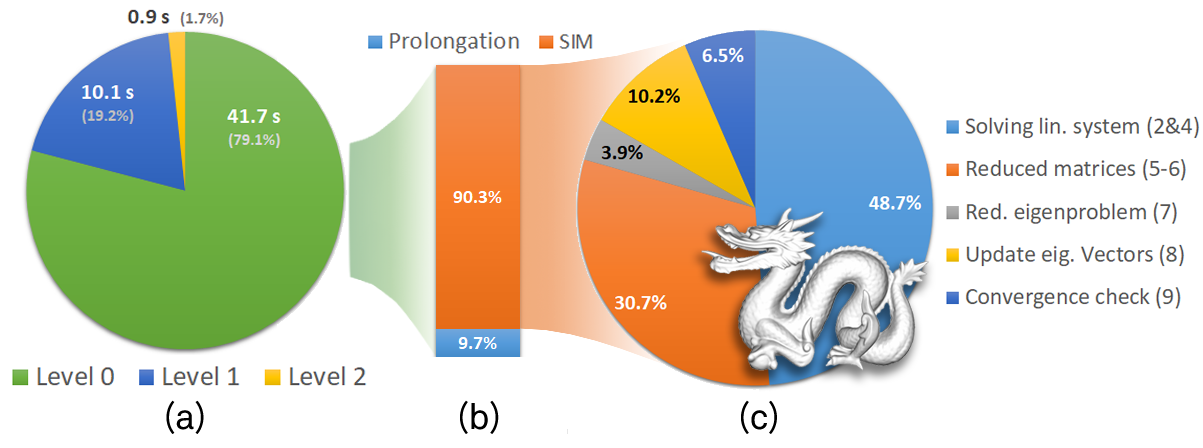}  
  \parbox[t]{1.0\columnwidth}{\relax  }
  \caption{\label{fig:timing}
  Analysis of timings of the HSIM for the Laplace--Beltrami eigenproblem. 
  Distribution of the runtimes to the three levels (a), split of the time spent at the finest level between the prolongation of the solution from level~1 and the subspace iterations at the finest level (b) and distribution of the time of the subspace iteration to the individual steps in Algorithm~\ref{alg:basic_SIM}. The 200 lowest eigenpairs are computed on the Dragon mesh with 150k vertices.  
}
\end{figure}

\paragraph{Timings} 
Table~\ref{tab:construction} lists timings of our HSIM implementation for the computation of the $p$ lowest 
eigenpairs of the discrete Laplace--Beltrami operator, eq. (\ref{eq.eigenProb}), for meshes with different 
sizes and values of $p$. 
Individual timings for hierarchy construction and solving the eigenproblems using the hierarchy are listed. Moreover, 
iteration counts for subspace iterations on the individual levels are provided. For the coarsest level, a dense 
solver is used instead of the subspace iteration, therefore, the table lists \emph{F}'s instead of a number for the 
coarsest level. The solver's convergence tolerance is set to $10^{-2}$ for all examples. 
In all cases, the required number of iterations on the finest level is reduced to one by the hierarchical approach. 
Figure~\ref{fig:timing} provides more details for one example, the computation of the lowest 200 eigenpairs on a dragon model with $150k$ vertices using a hierarchy with three levels. 
The figure shows (a) how the runtimes split over the different levels of the hierarchy, (b) for the finest level the division between prolongation of the solution for the second finest level and subspace iterations, and (c) the breakdown of the timings of the individual steps of the subspace iterations (Algorithm~\ref{alg:basic_SIM}) on the finest level. 
The figure illustrates that, when three levels are used, most of the runtime is spent on the finest level, almost 80\% for the shown example, and that the restrictions of the stiffness and mass matrices and solving the linear systems are the most costly steps of HSIM.  

\input{tables/MLsolver.tex}

\begin{figure}[t]
  \centering
  \includegraphics[width=1.0\linewidth]{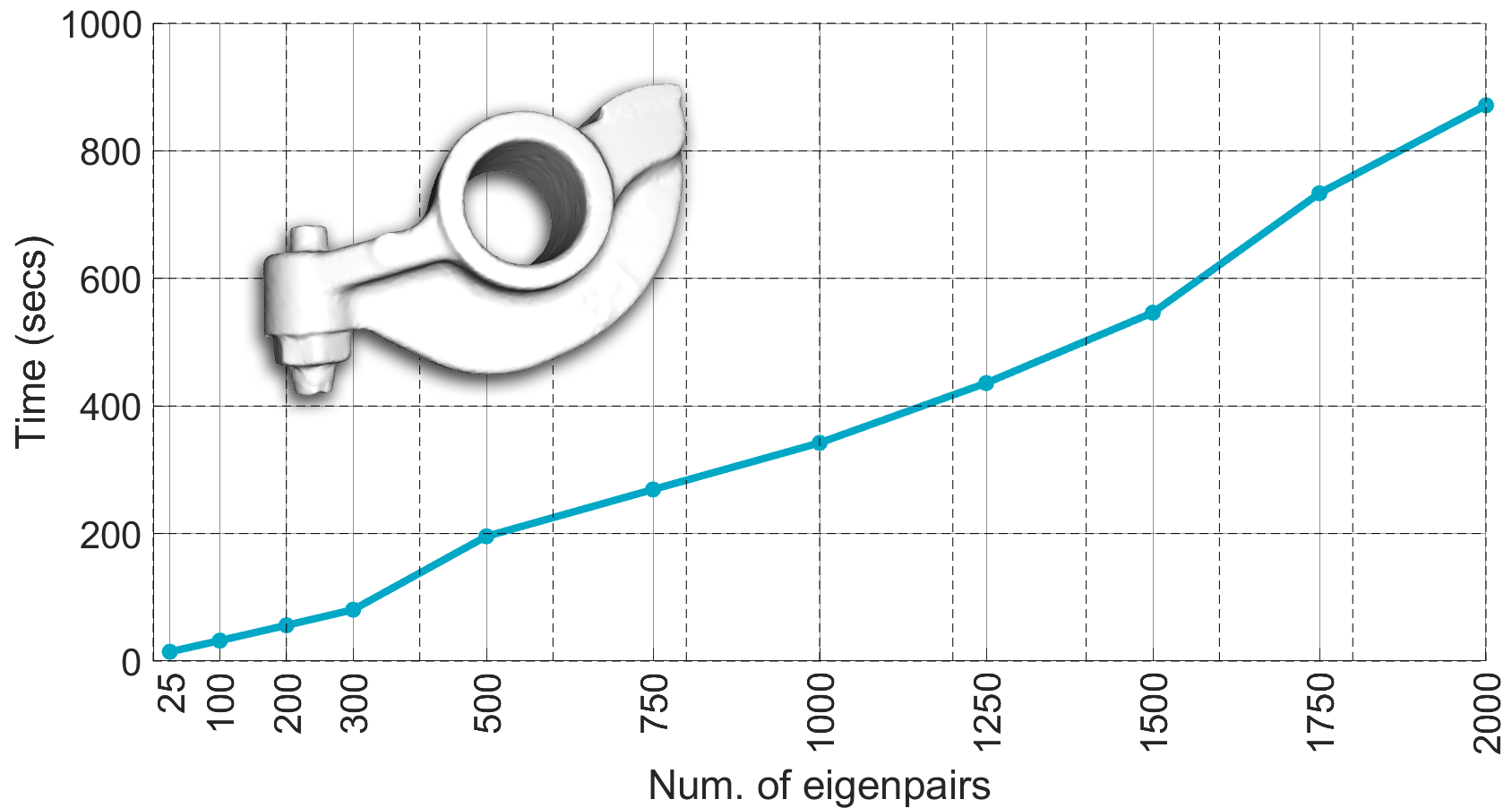}
  \parbox[t]{1.0\columnwidth}{\relax  }
  \caption{Plot listing the runtime of HSIM over the number of Laplace--Beltrami eigenpairs to be computed on the Rocker Arm model with 270k vertices.
  }
\label{fig:eigenvalues}
\end{figure}

\begin{figure}[b]
  \centering
  \includegraphics[width=1.0\linewidth]{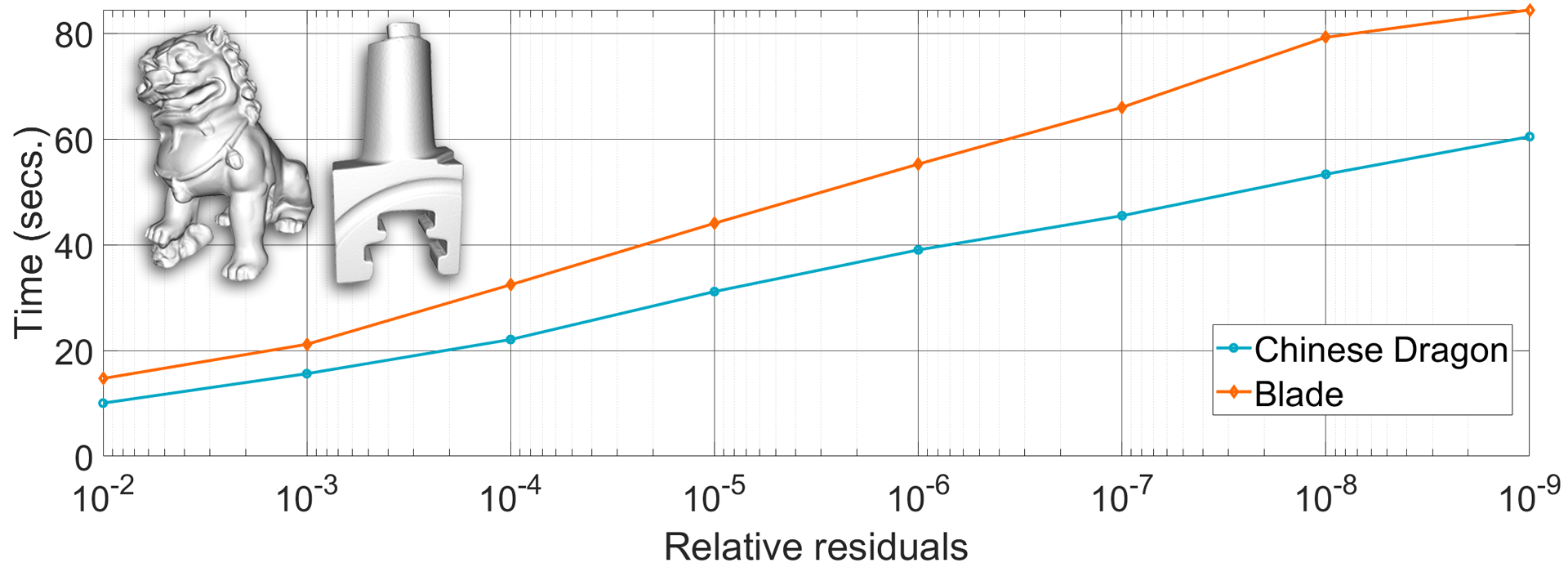}
  \parbox[t]{1.0\columnwidth}{\relax  }
  \caption{\label{fig:residue}
	The plot of the required computation time for different error tolerances when computing the first 100 eigenpairs 
	of Laplace--Beltrami operator on the Chinese Dragon (127k vertices) and on the Blade model (200k vertices).
}
\end{figure}

Figure~\ref{fig:eigenvalues} lists runtimes for different numbers of eigenpairs to be computed. 
In our experiments, we found that the runtime grows linearly even when computing several thousand eigenpairs. This is illustrated by the timings listed in the figure. We expect this linear trend to continue as long as the runtime is dominated by the time needed for the solving of the linear systems (step 4 of Algorithm~\ref{alg:basic_SIM}). 
At some point, solving the dense eigenproblems (step 7 of Algorithm~\ref{alg:basic_SIM}), which does 
not scale linearly with the number of eigenpairs, will be the most expensive step and the trend will no longer be linear.

\input{tables/vsAll}

For most experiments, we set the convergence tolerance, $\epsilon$ in Algorithm~\ref{alg:ours}, to $10^{-2}$. 
Figure~\ref{fig:residue} lists runtimes over the convergence tolerance for the computation of 100 eigenpairs on two different meshes, the Blade model with 200k vertices and the Chinese Dragon with 127k vertices.
The figure illustrates that low tolerances such as $10^{-9}$ can be achieved and that the time grows proportional with the relative residual. 
Roughly speaking, we observe in our experiments that the number of iterations that are needed on the finest grid grows by two for a decrease of one order of magnitude in the relative residual.

\begin{figure}[b]
    \subfloat[\minorrev{HSIM, 10k}\label{fig:HSIM_10k}]{
        \includegraphics[width=0.49\linewidth]{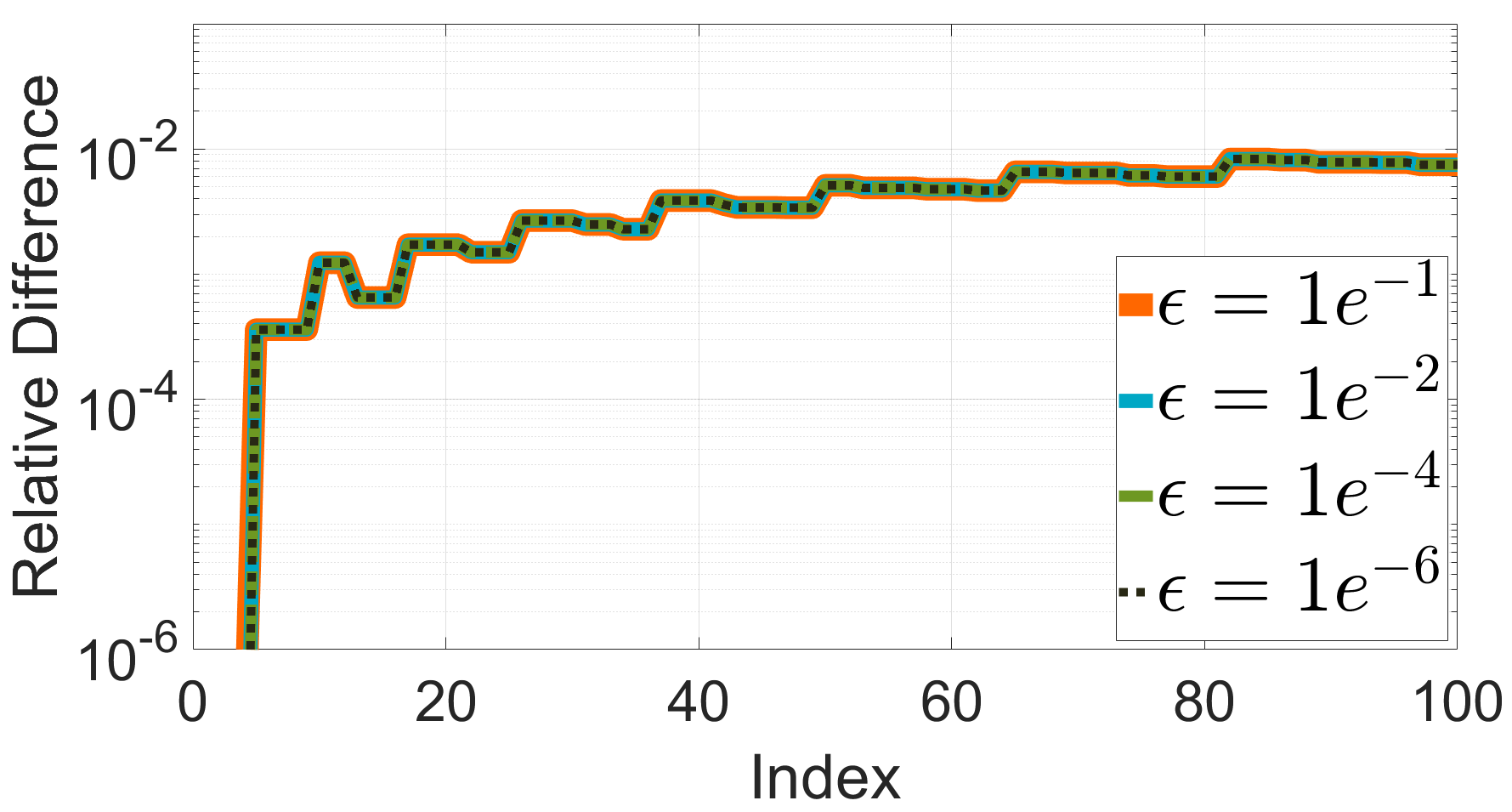}
    }
    \subfloat[\minorrev{HSIM, 100k} \label{fig:HSIM_100k}]{
        \includegraphics[width=0.49\linewidth]{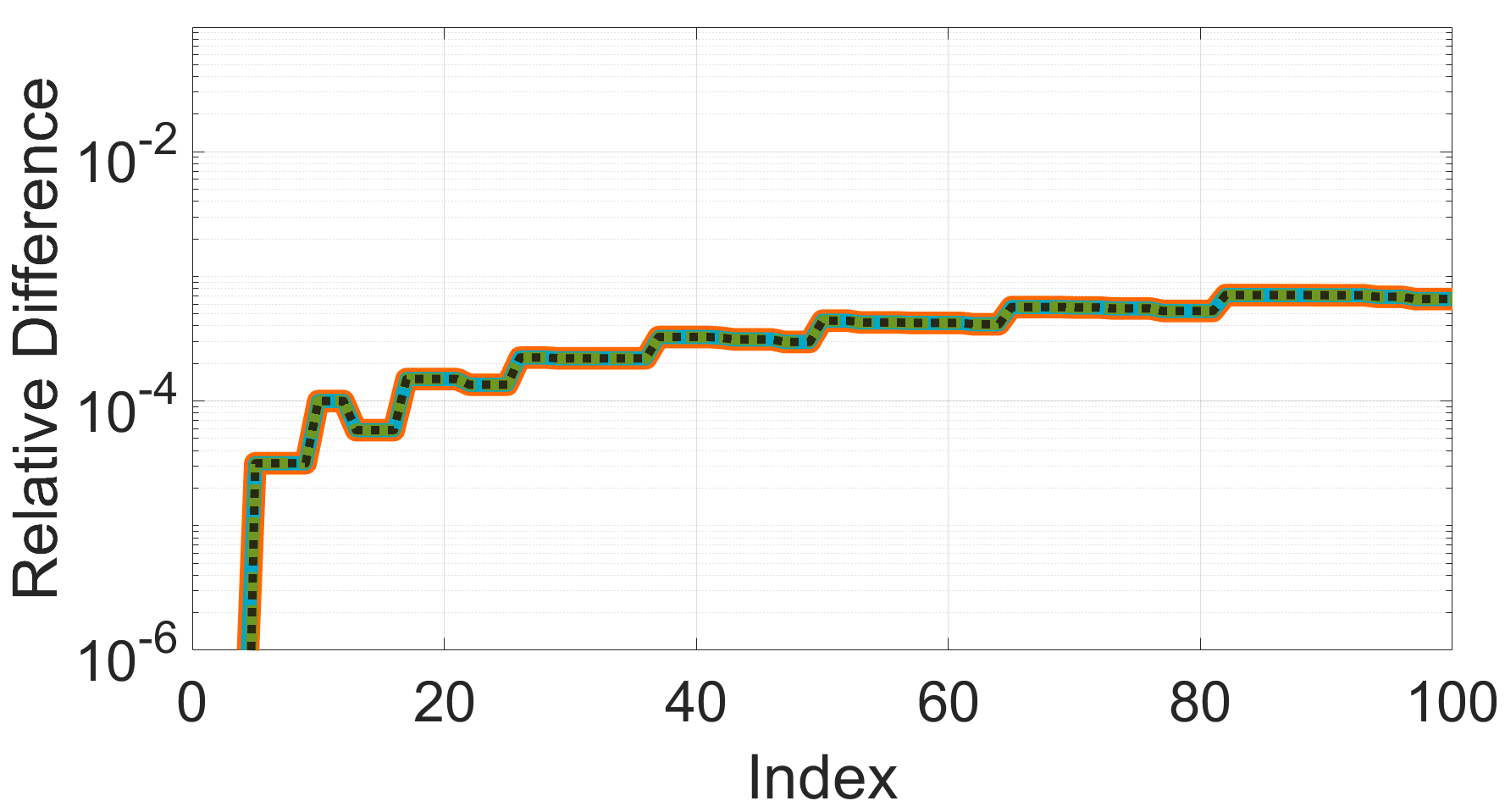}
    }
    \\
    \subfloat[\minorrev{HSIM, 1m} \label{fig:HSIM_1m}]{
        \includegraphics[width=0.49\linewidth]{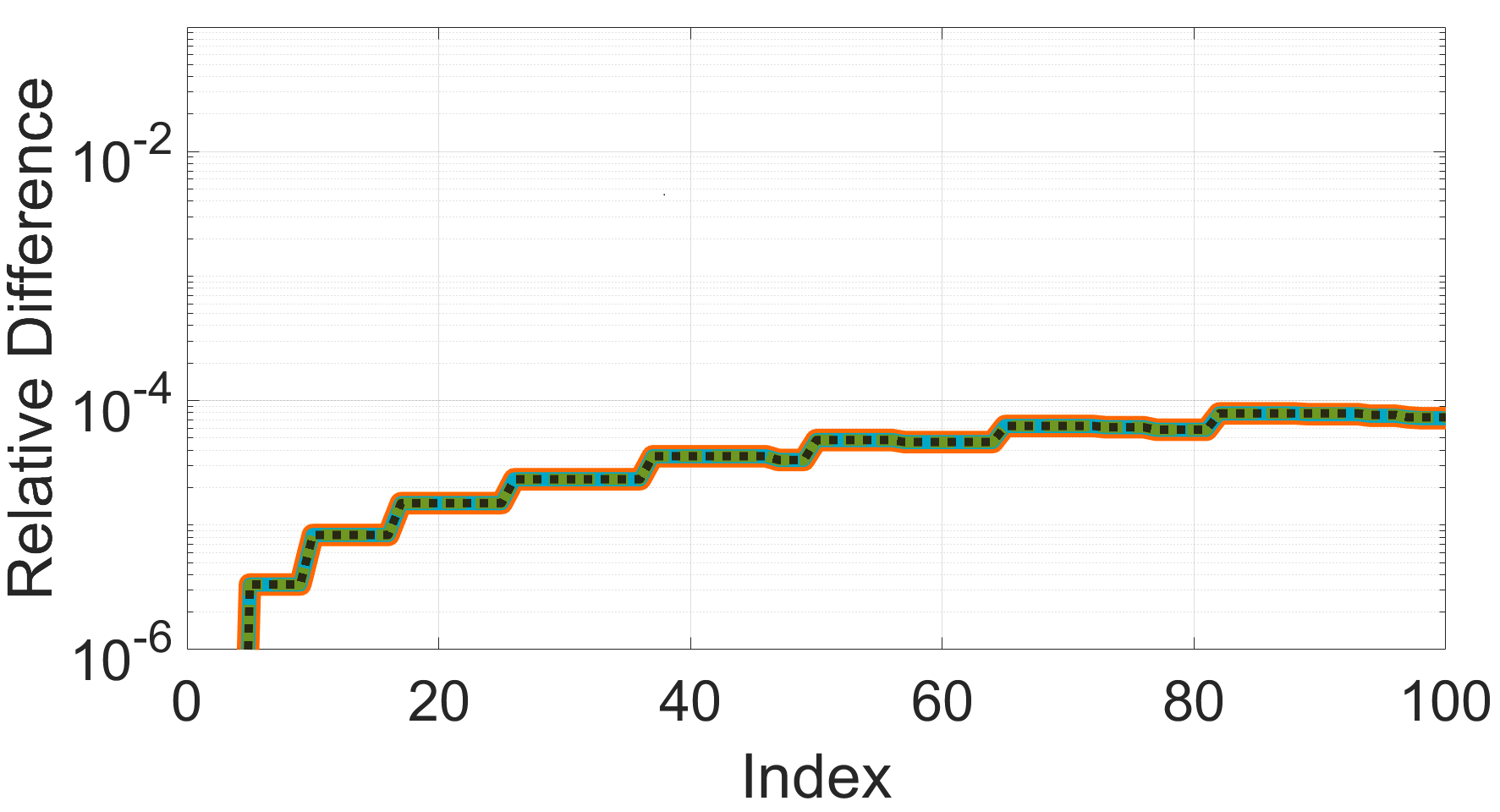}
    }
    \subfloat[\minorrev{IPM, 100k} \label{fig:IPM_100k}]{
        \includegraphics[width=0.49\linewidth]{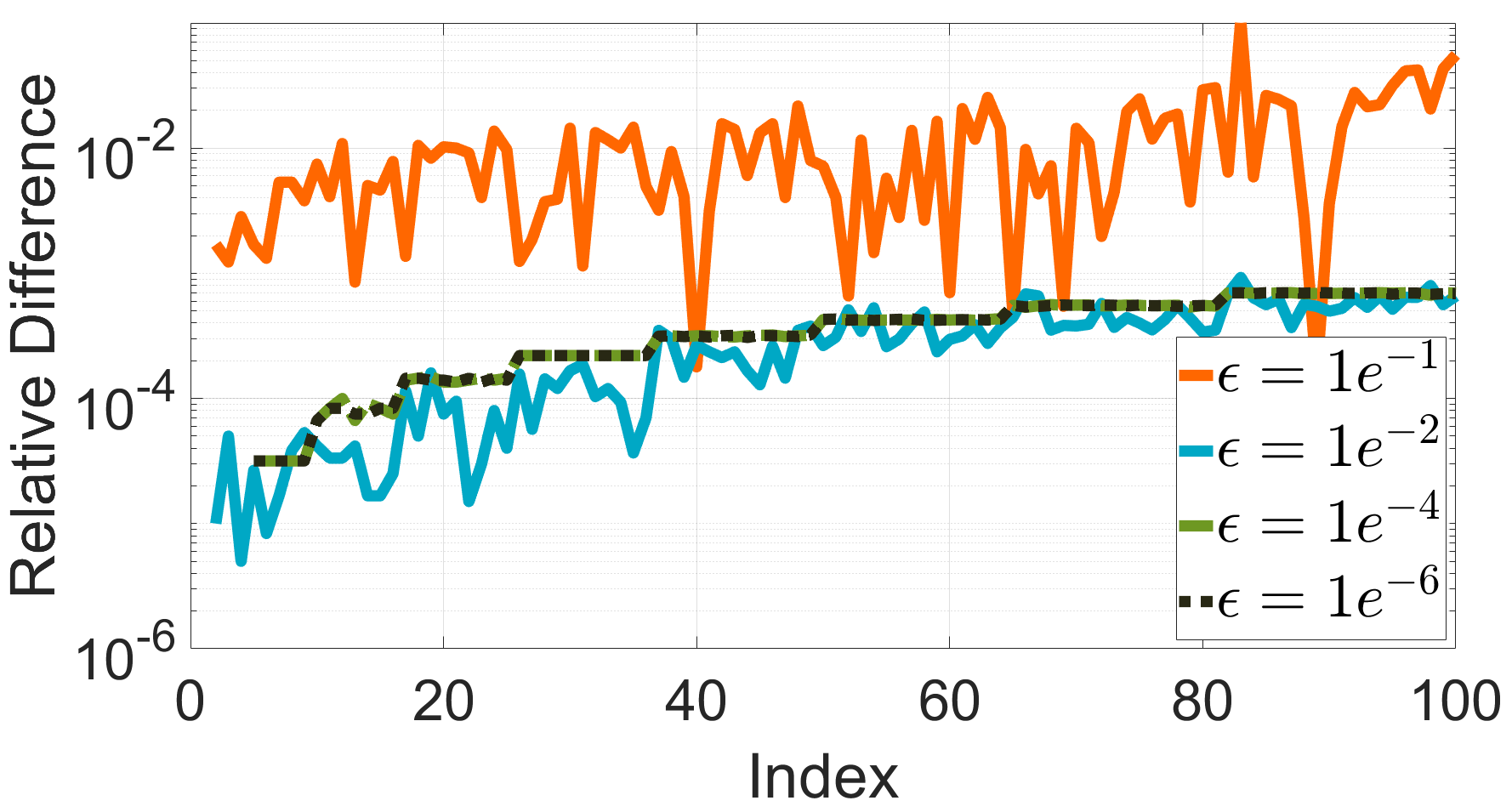}
    }
  \caption{\label{fig:error}
  \revised{Relative difference of numerical approximations of the eigenvalues of the unit sphere to the analytic solutions are shown.}
}
\end{figure}

\begin{figure}[b]
    \subfloat[\minorrev{Eigenvalues}\label{fig:edge-flip-1}]{
        \includegraphics[width=0.49\linewidth]{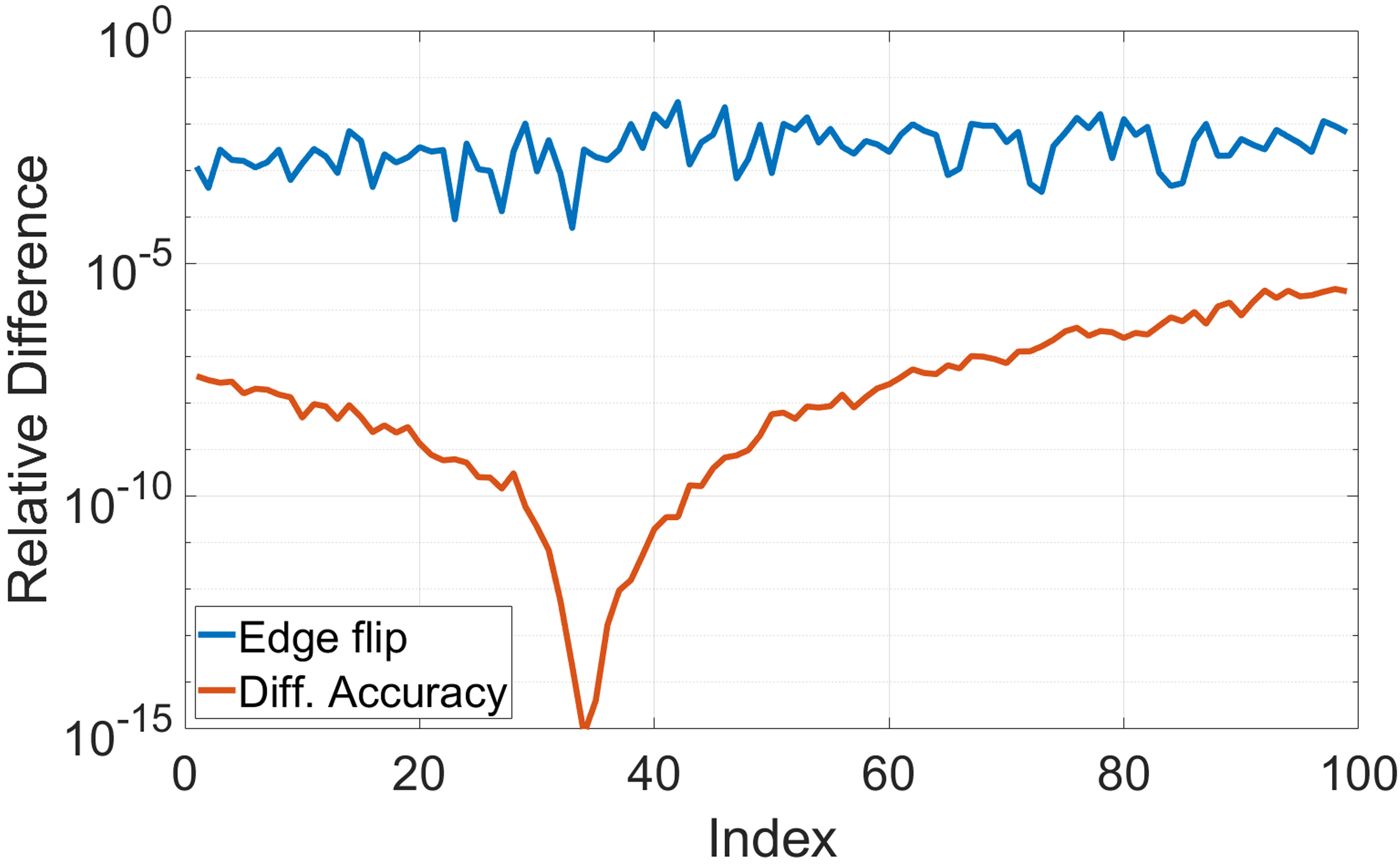}
    }
    \subfloat[\minorrev{Eigenvectors}\label{fig:edge-flip-2}]{
        \includegraphics[width=0.49\linewidth]{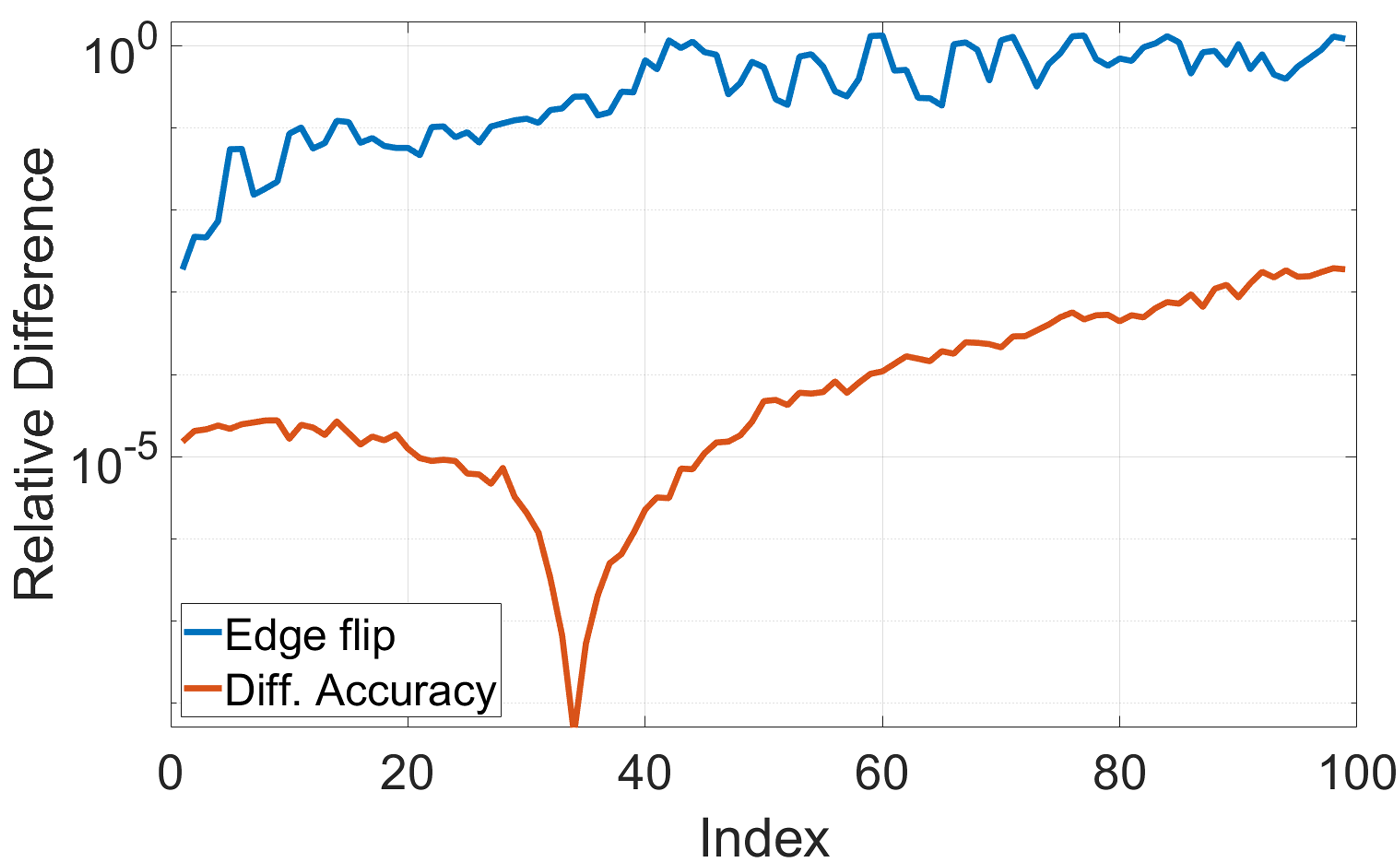}
    }
    \\
    \subfloat[\minorrev{Model (original)}\label{fig:edge-flip-3}]{
        \includegraphics[width=0.49\linewidth]{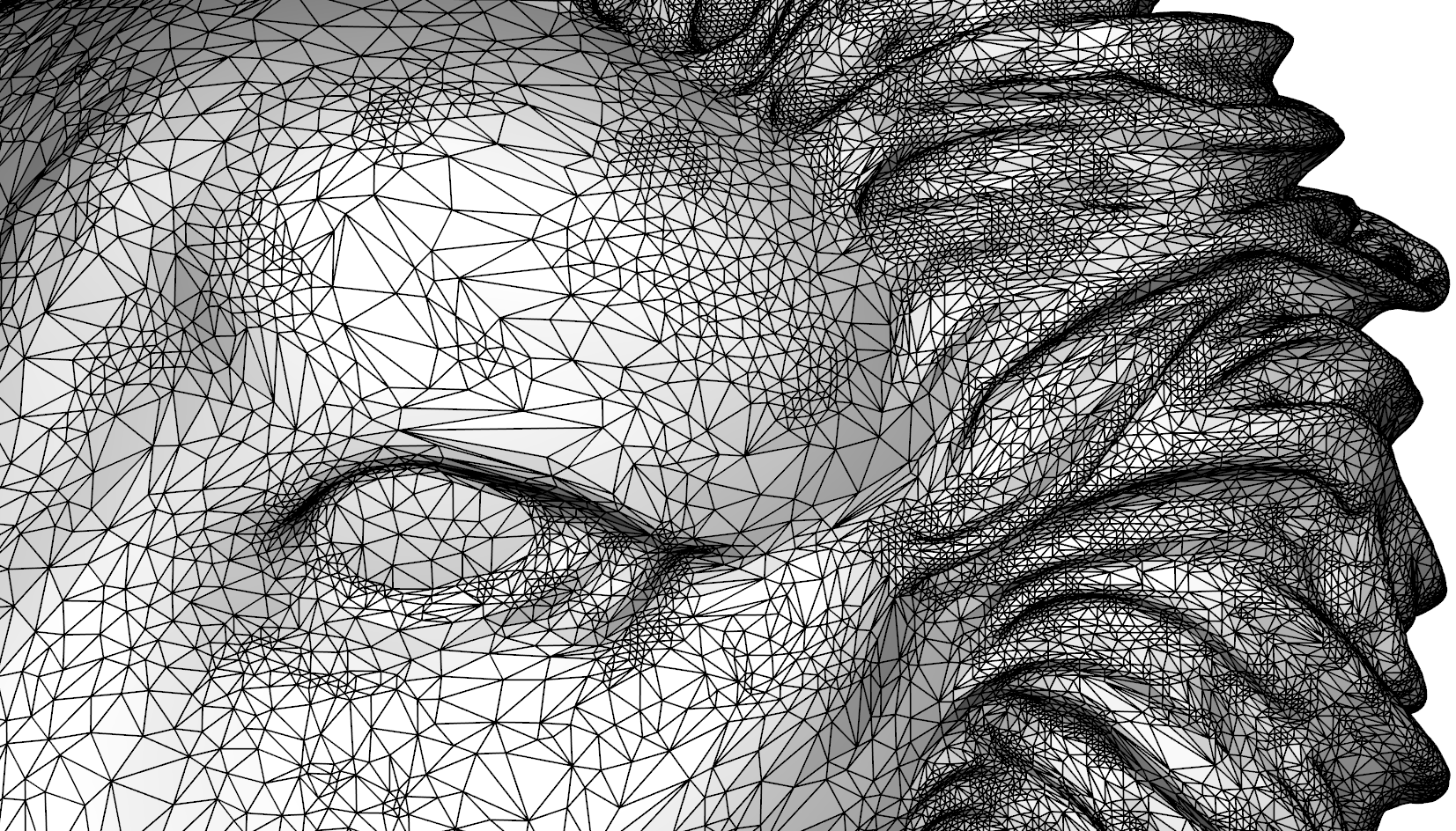}
    }
    \subfloat[\minorrev{Model (edge flips)}\label{fig:edge-flip-4}]{
        \includegraphics[width=0.49\linewidth]{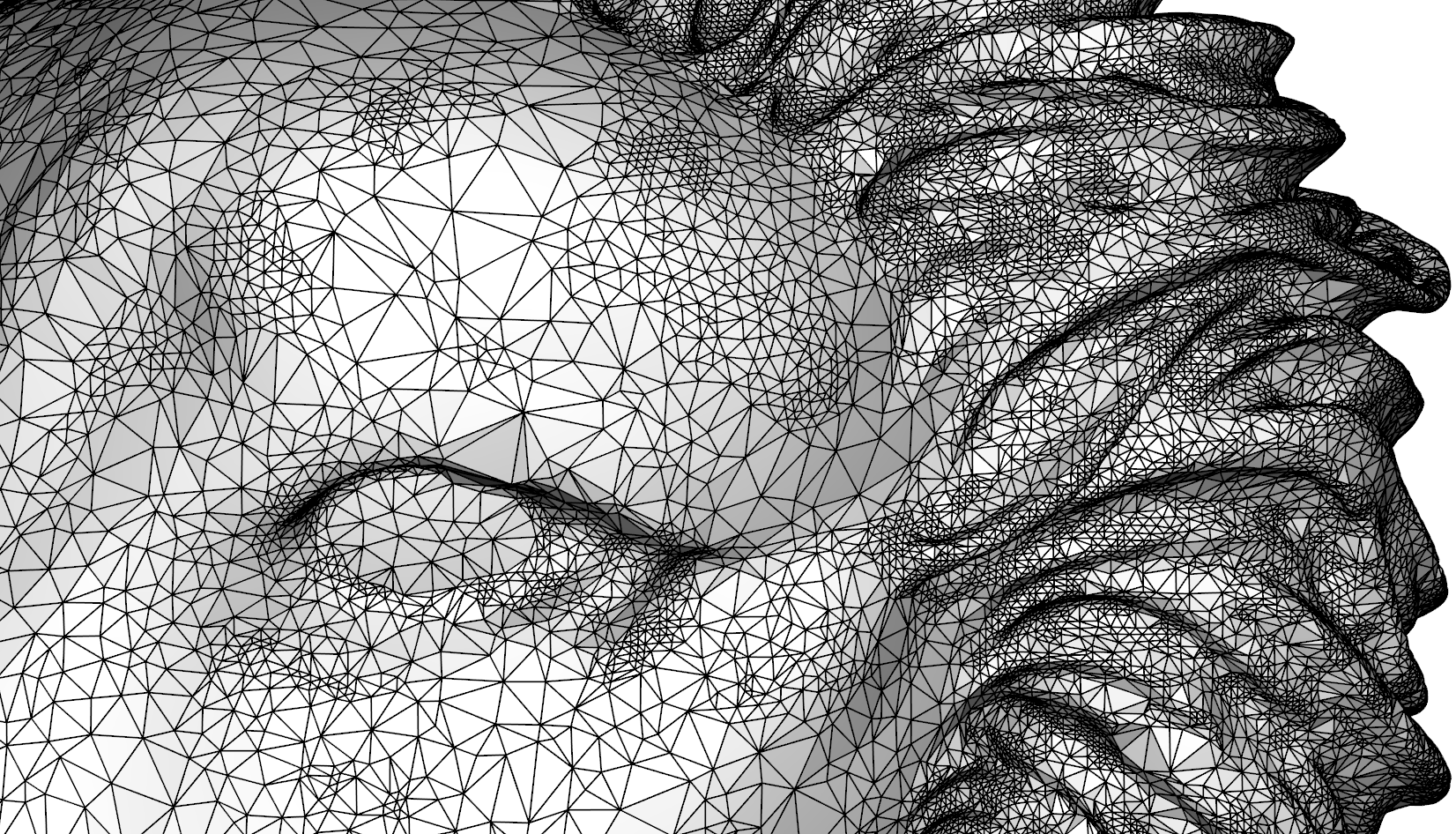}
    }
     \caption{\revised{Comparison of the relative difference of the eigenvalues and the eigenvectors between two meshes that approximate the same surface (blue graph) and solutions for different convergence tolerance on one of the meshes (red graph).}}
     \label{fig:edge-flip} 
\end{figure}

\revised{
\paragraph{Termination criterion}
To test for convergence, see line 9 of Algorithm 1, we use the criterion stated in (\ref{eq.conv2}). This test ensures the convergence of the eigenvalues as well as the convergence of the eigenvectors. 
To determine a suitable value for the convergence tolerance $\varepsilon$, we performed several experiments. 
We discuss two experiments in this paragraph, the supplementary material includes additional experiments.
Based on the results of these experiments, we used a tolerance of $\varepsilon=10^{-2}$ for the evaluation of HSIM. 
In the first experiment, we consider three different discretizations of the unit sphere with regular meshes (having 10k, 100k and 1m vertices) and measure the difference between the computed eigenvalues for different tolerances ($\varepsilon=10^{-1},10^{-2},10^{-4},10^{-6}$) and the analytical solution. 
The results are shown in Figure~\ref{fig:error}.
For all three discretizations, the difference between the numerical solutions for different tolerances is small compared to the approximation error, that is, the difference to the analytical solution. 
We would like to note that the convergence test establishes an upper bound on the convergence of the eigenpairs. In particular, for the lowest tolerance, $\varepsilon=10^{-1}$, the solutions computed by HSIM are often already more accurate when the process terminates. 
One reason for this is that the method terminates only after all eigenpairs pass the convergence test. 
We therefore conducted an additional experiment using the inverse power method to compute the eigenpairs one by one and stop the iteration for each eigenpair when the convergence tolerance is reached. The results are shown in Figure~\ref{fig:error}~(d). In this experiment, differences in accuracy occur between the numerical solution for $\varepsilon=10^{-1}$ and the other solutions ($\varepsilon=10^{-2},10^{-4},10^{-6}$), which indicates that a tolerance of $\varepsilon=10^{-1}$ is not sufficient.

In a second experiment, we compute eigenpairs for two different meshes approximating the same surface. The second mesh was created by flipping edges of the first mesh. 
For both meshes, we compute the lowest eigenpairs for the tolerance $\varepsilon=10^{-2}$ and as reference for $\varepsilon=10^{-8}$.
Since the two meshes have the same vertices, we can compare both the eigenvalues and the eigenvectors. 
Figure~\ref{fig:edge-flip} shows the difference between the reference solutions ($\varepsilon=10^{-8}$) on both meshes (blue graph) and for one mesh, the difference between the solutions for $\varepsilon=10^{-2}$ and $\varepsilon=10^{-8}$ (red graph). }
\minorrev{
To measure the difference of eigenvectors the relative $L^2$-norm is used. 
}
\revised{
It can be seen that the difference between the reference solutions on the two meshes is more than three orders of magnitude larger than the difference between the solutions for different tolerances.

We want to note that the convergence test (\ref{eq.conv2}) does not directly measure the deviation from the exact solution. In our experiments (for example in Figure~\ref{fig:edge-flip}), we see that the relative difference between the solution for a tolerance of $\varepsilon=10^{-2}$ and the reference solution, which is computed with $\varepsilon=10^{-8}$, is usually much smaller than $10^{-2}$. 
In Figure~\ref{fig:edge-flip}, and also in Figure~\ref{fig:vsBMG}, one can observe that the errors generated by HSIM are smaller for the eigenvalue pairs whose index is about one-third of the total number of computed eigenvalues than for the others. 
This is due to our shifting strategy, which makes these eigenpairs converge faster.
}

\input{tables/levels}

\input{tables/support.tex}

\input{tables/shift.tex}

\minorrev{When evaluating the convergence criterion, equation (\ref{eq.conv2}), the standard norm of $\mathbb{R}^n$ is commonly used to replace the $M^{-1}$-norm. 
The reason is that the evaluation of the $M^{-1}$ norm can be costly. 
For our experiments, we used the $M^{-1}$-norm at the finest level as the $M$ matrices are diagonal, and, therefore, can be easily inverted. 
At the coarser levels, however, the restricted 
matrices $M^{\tau}$ (see line 6 of Algorithm~\ref{alg:ours}) are no longer easy to invert. To save the effort of computing a factorization of the $M^{\tau}$ matrices, we replace the $M^{-1}$-norm by the standard norm for the convergence check on all but the finest level. We want to emphasize that since on the finest level we use the $M^{-1}$-norm, the error tolerances are respected. The simplification could only lead to more or fewer iterations on the coarser levels. Since the convergence check uses the relative norm, a global scaling factor to better match the standard norm and the $M^{-1}$-norm is not needed. 
We did not observe differences in the numbers of iterations, when using the standard norm instead of the $M^{-1}$-norm for the convergence test on the coarser levels in our experiments.}

\paragraph{Number of levels} \label{par:level}
A parameter HSIM needs as user input is the number of levels of the hierarchy, see Algorithm~\ref{alg:ours}. 
By increasing the number of levels, one can reduce the number of iterations required on the finest grid. 
On the other hand, increasing the number of levels leads to additional computational costs on the levels below the finest level. 
Table~\ref{tab:level} lists computation times and iteration counts for the individual levels for computations with different meshes sizes, number of eigenpairs and convergence tolerances. 
For most of these examples, three levels yield the shortest runtime.

\paragraph{Support region}
For the construction of the prolongation matrices~$U^{\tau}$, the radius of its domain of influence, $\rho^{\tau}$, must be defined 
individually for each level.
We use eq.~(\ref{eq.radius}), which allows us to set the radii on all levels by means of a control parameter $\sigma$. 
\setlength{\columnsep}{5pt}
\setlength{\intextsep}{0pt}
\begin{wrapfigure}{r}{0.33\linewidth}	
    \centering
    \includegraphics[width=1.0\linewidth]{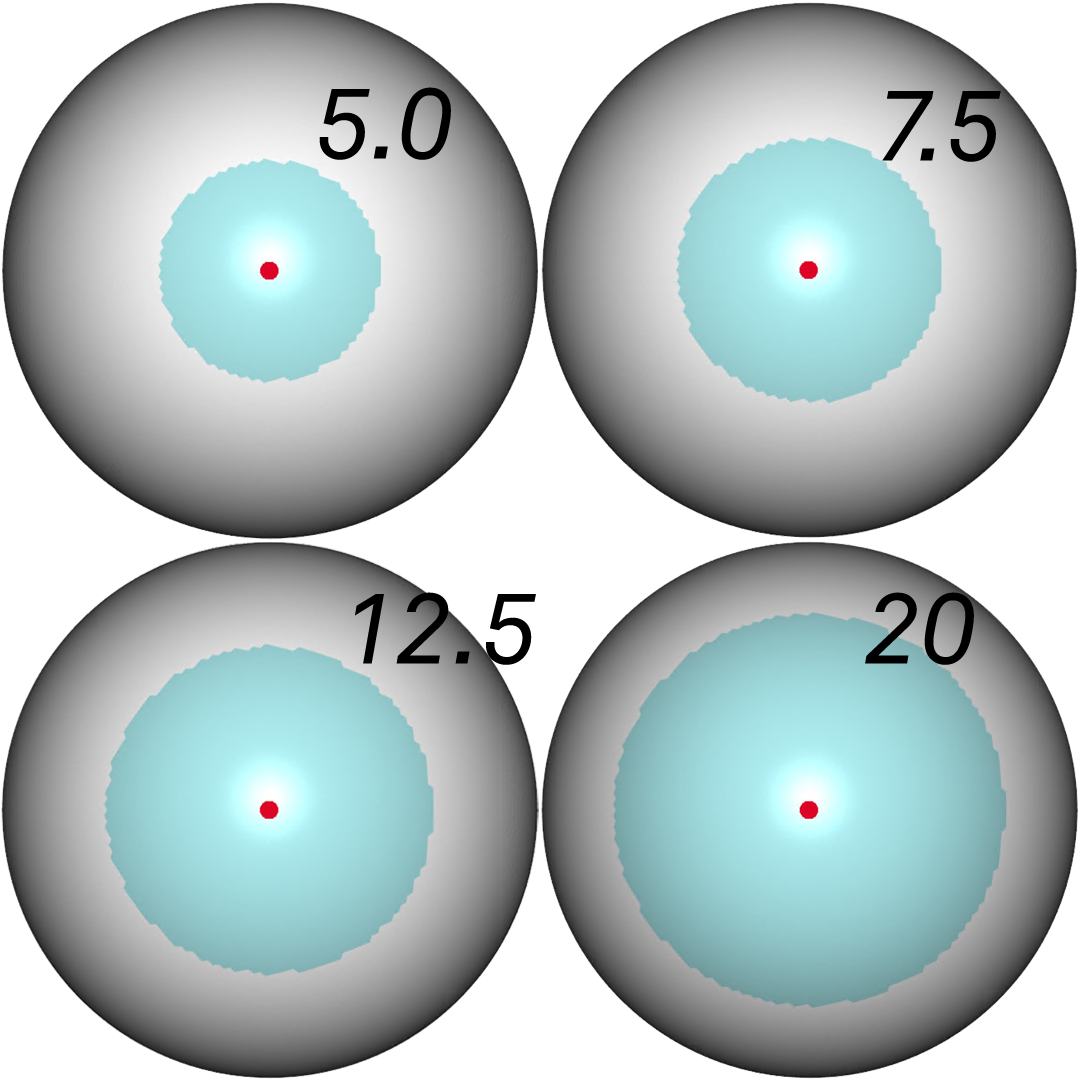}
	\vspace{-10pt}
	\label{fig:var-support}
\end{wrapfigure}
This value is the average expected number of non-zero entries per row of the matrices $U^\tau$. The inset figure shows the areas of influence around one point for different values $\sigma$.  
A smaller value for $\sigma$ results in matrices $U^{\tau}$ with less non-zero entries and thus less computational effort per iteration. On the other hand, a too small value for $\sigma$ can increase the number of iterations needed on each level. 

In our experiments, we have identified a value of $\sigma=7$ as a good trade-off. This means that in each level, each vertex of $V^{\tau}$ in average is coupled to six neighbor vertices, which agrees with the average valence in a triangle mesh. 
Table~\ref{tab:support} shows iteration counts and runtimes for different values of $\sigma$ for eigenproblems on two meshes with 200k and 475k vertices and different numbers of eigenpairs to be computed. 
The value $\sigma=7$ reaches in all cases either the lowest runtime or a time close to the lowest runtime.

\input{tables/boundary}


\begin{figure*}[t]
  \subfloat[Sphere, 160k] {
    \includegraphics[width=0.475\linewidth]{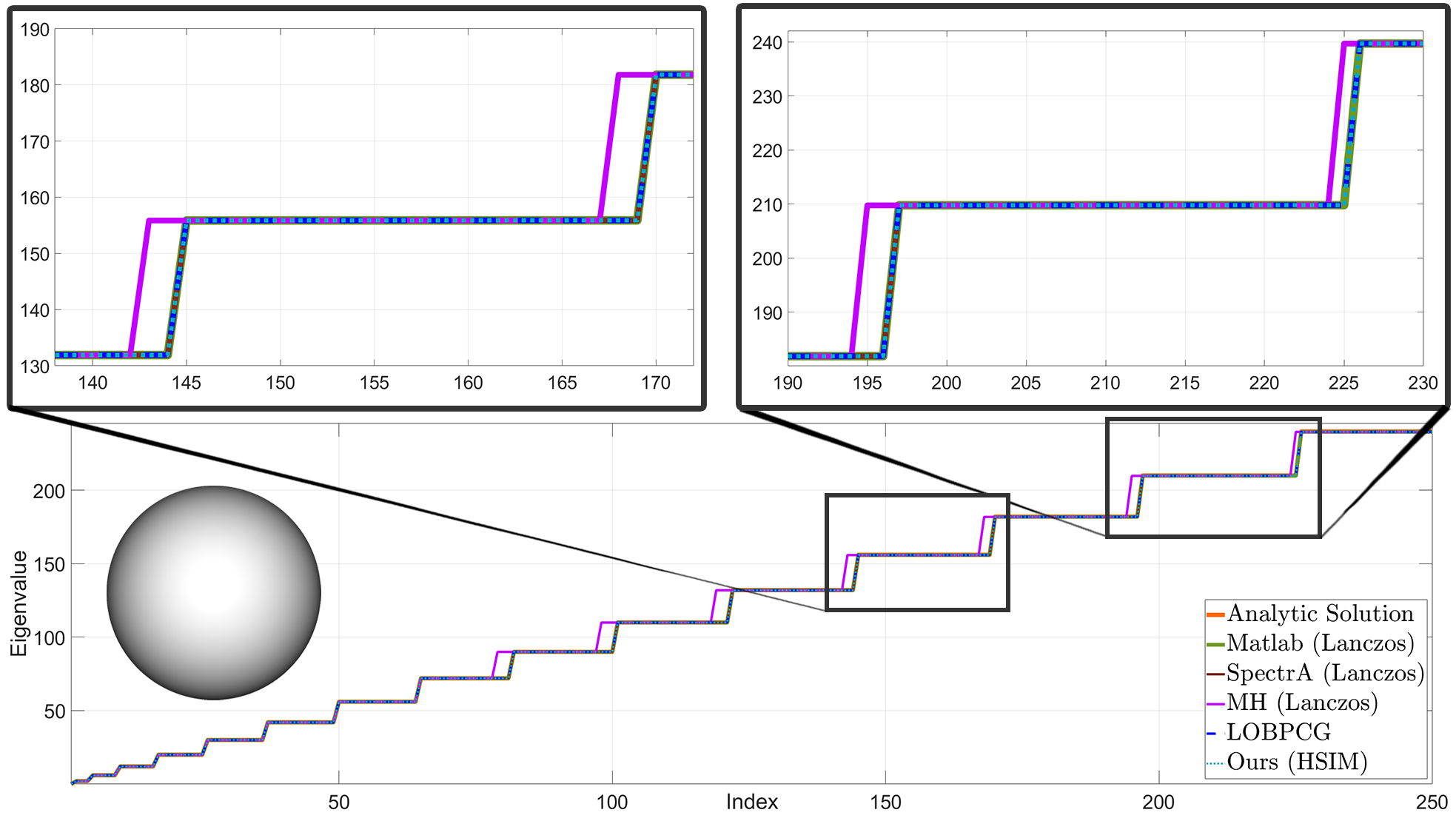}\label{fig:eigvals-1}
  }
  \hfill
  \subfloat[Ball, 90k]{
    \includegraphics[width=0.475\linewidth]{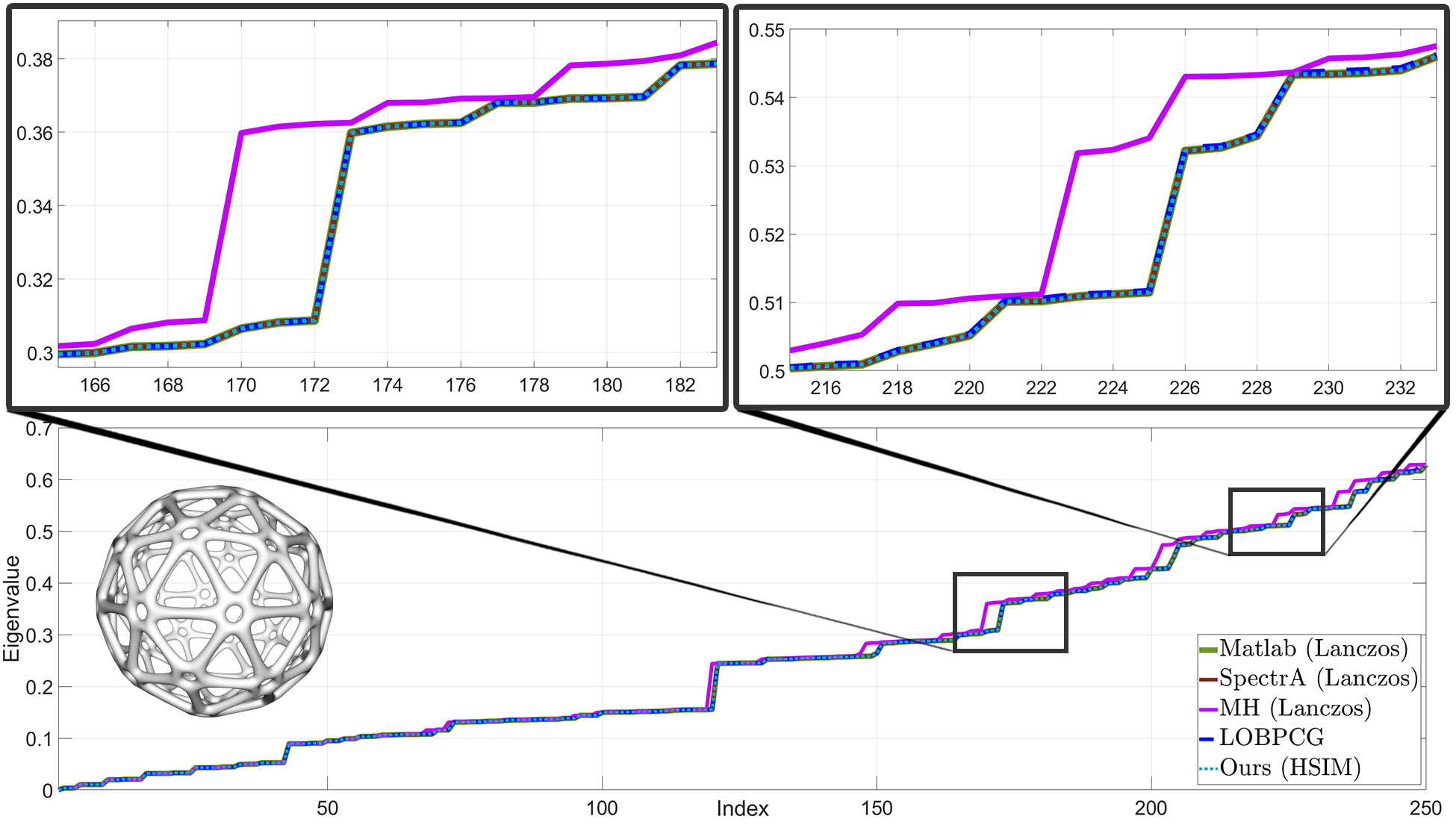}\label{fig:eigvals-2}
  } 
  \parbox[t]{1.0\columnwidth}{\relax  }
  \caption{\revised{
	The lowest 250 Laplace--Beltrami eigenvalues computed with HSIM and three different Lanczos solvers on a discrete sphere with 160k vertices (left) and a surface with many symmetries and 90k vertices (right). For the sphere, the analytic solution is shown as a reference.}} 
  \label{fig:lanczos}
\end{figure*}

\begin{figure}[b]
    \subfloat[Dirichlet\label{fig:boundary-dirichlet}]{
            \includegraphics[width=0.49\linewidth]{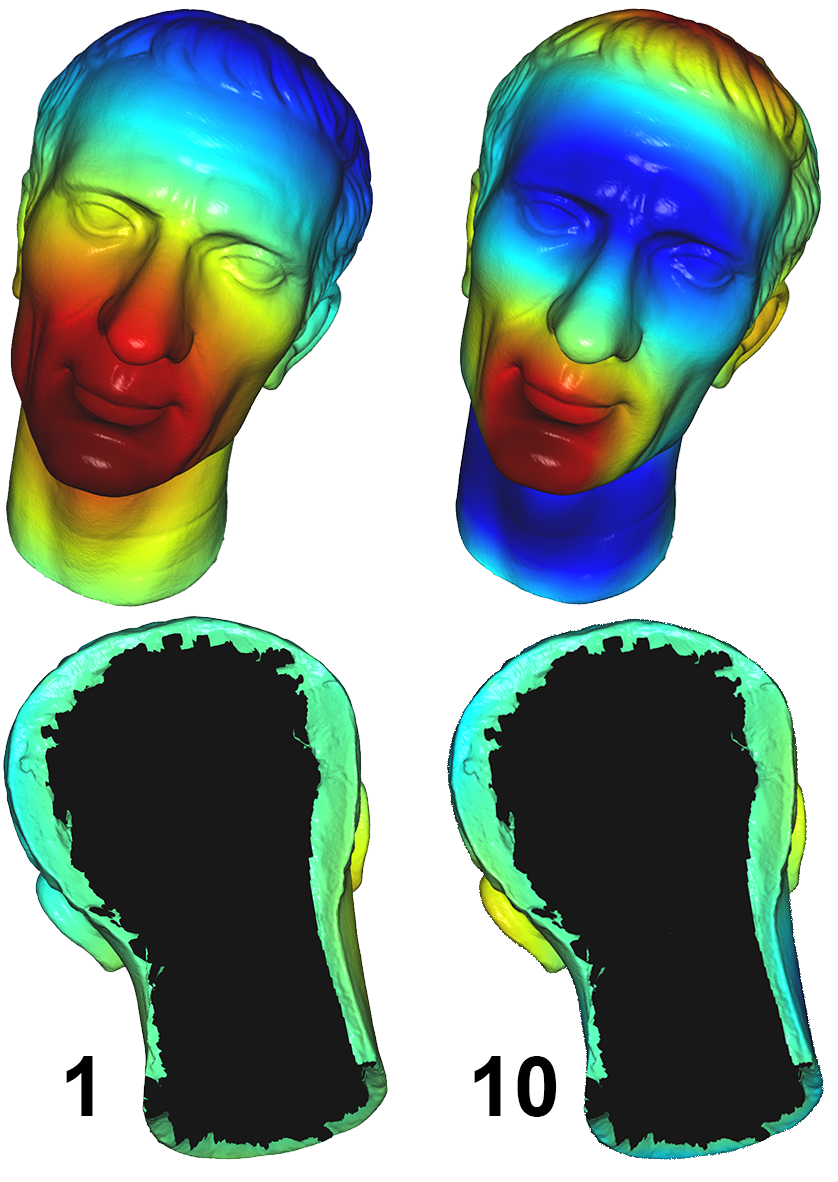}
    }
    \subfloat[Neumann \label{fig:boundary-neumann}]{
        
            \includegraphics[width=0.49\linewidth]{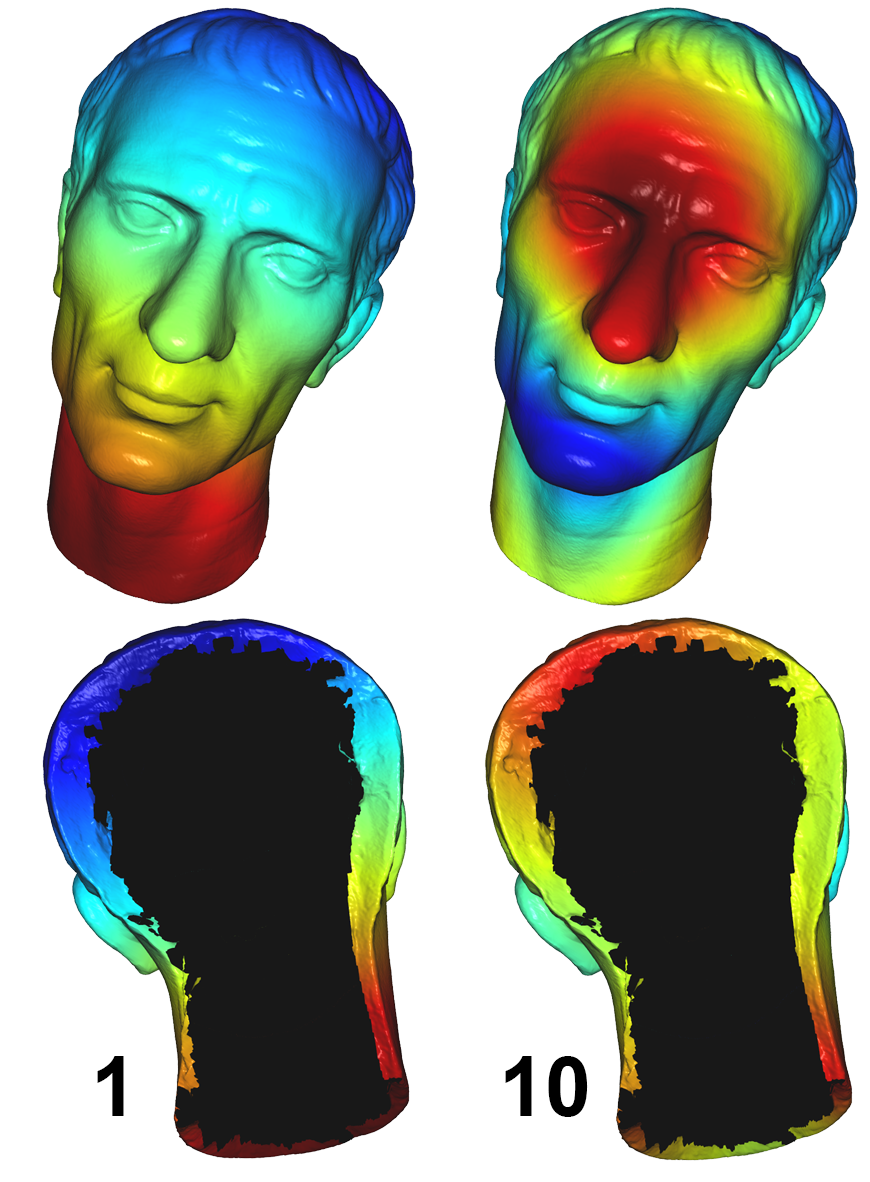}
    }
     \caption{The first (left) and tenth (right) eigenfunction of the Laplace--Beltrami operator on a surface with boundary using Dirichlet and Neumann boundary conditions are shown.}
     \label{fig:boundary} 
\end{figure}

\paragraph{Shifting strategy}
Matrix shifting can reduce the number of required subspace iterations on all levels. 
In our experiments, we use a heuristic, which is described in Step 12 of Algorithm~\ref{alg:ours}, to automatically determine the shifting parameter $\mu$.
This heuristic is based on the aggressive shifting technique from \cite{zhao2007accelerated}. 
We set $\mu$ equal to the current approximate eigenvalue $\Lambda_{jj}$ with index $j={\lfloor \alpha p \rfloor}$. 
Here, $\alpha$ should take a value between 0 and $0.5$. In
Algorithm~\ref{alg:ours}, we set $\alpha=0.1$. 
Table~\ref{tab:shift} lists iteration counts for different values of $\alpha$. Results for different numbers of eigenpairs and error margins are shown.
We used values between $0.1$ and $1/3$ for $\alpha$ in our experiments. 

\revised{
\paragraph{Surface with boundary}
We applied HSIM to the computation of Laplace--Beltrami eigenproblems on surfaces with boundary. 
We experimented with Dirichlet and Neumann boundary conditions and used the same hierarchy and basis construction as for surfaces without boundary. 
Examples of eigenfunctions on surfaces with boundary are shown in Figure~\ref{fig:boundary}. 
Table~\ref{tab:boundary} shows for an example mesh the runtimes and iteration counts for Dirichlet and Neumann boundary conditions. 
The runtimes are comparable to the runtimes we observe for meshes without boundary and a comparable number of vertices.
}

\section{Comparisons}\label{sect.comparison}
In this section, we discuss comparisons of HSIM to alternative methods. 
Laplace--Beltrami eigenproblems are commonly solved in graphics applications 
using Lanczos methods \cite{Vallet2008}. Therefore, we begin this section with the 
comparison to Lanczos solvers. An alternative to Lanczos schemes is the SIM \cite{bathe2013subspace}. Since 
HSIM is based on SIM, this comparison provides a basis to quantify the gains 
resulting from our hierarchy. \revised{The third solver to which we compare HSIM is the 
Locally Optimal Block Preconditioned Conjugate Gradient Method~\cite{knyazev2001toward}.}
Lastly, we compare HSIM to the multilevel correction scheme that was introduced in~\cite{Chen2016,Lin2015}.
\revised{
In our comparisons, we use the same convergence test (\ref{eq.conv2}) for all methods and set the tolerance to $\varepsilon=10^{-2}$, which is the value we determined in our experiments, see Section~\ref{sect.experiments}.
}\minorrev{To implement this for the methods we compare to, we check for
convergence after every iteration and stop when the convergence test is passed.
Once the required number of iterations is known, we re-run the computation
without convergence test and record the timings.
}

\paragraph{Lanczos schemes}
Schemes based on Lanczos iterations are commonly used for solving large-scale, 
sparse, symmetric eigenproblems. These methods have been studied and 
improved over decades. \textsc{Arpack}'s implementation of 
the Implicitly Restarted Lanczos Method is well-established \cite{Lehoucq1998}. We compare HSIM 
with \textsc{Matlab}'s \texttt{eigs} \minorrev{(\textsc{Matlab}~9.8, R2020a)} that interfaces \textsc{Arpack} and with 
\textsc{SpectrA} \cite{Qiu2015} that offers an alternative implementation of 
the Implicitly Restarted Lanczos Method. In addition to that, we compare to 
the authors' implementation of the band-by-band, shift-and-invert Lanczos 
solver that was introduced in \cite{Vallet2008}. We denote this solver by Manifold Harmonics (MH).
If a diagonal, or lumped, mass matrix is used in (\ref{eq.eigenProb}), the generalized 
eigenproblem can easily be transformed to an `ordinary' eigenproblem as described in \cite{Vallet2008}. 
We have tested all three Lanczos solvers on the `ordinary' eigenproblem. 

The runtimes for meshes of different complexity and with different numbers of 
eigenpairs are given in Table~\ref{tab:vsAll}. The listed runtimes for HSIM 
also include the construction of the hierarchy and prolongation operators. In 
our experiments, HSIM was consistently faster than all three Lanczos schemes. 
This is also reflected in the table where HSIM is the fastest method for all 
combinations of mesh complexity and numbers of eigenpairs. In 
particular, for the difficult cases where a larger number of 
eigenpairs is computed, HSIM is significantly faster. 

Figure~\ref{fig:lanczos} shows plots of the lowest eigenvalues for two surfaces computed with different solvers. 
On the left side of the figure, 
numerical approximations of the eigenvalues of the unit sphere computed with 
the different solvers on a mesh with 320k triangles approximating the sphere 
are shown. For reference, the analytical solution is included to the plot. 
On the right side of the figure, results for a surface that exhibits different symmetries are shown. 
\textsc{SpectrA} and \textsc{Matlab} applied to the ordinary eigenproblem provided accurate results in our experiments that for the sphere example well-approximate the analytic solution.  
The results obtained with HSIM match the accuracy of \textsc{SpectrA} and \textsc{Matlab}. 
The band-by-band, shift-and-invert solver \cite{Vallet2008} meets the convergence tolerance for the 
individual eigenpairs, but some eigenpairs are skipped. This seems to happen at 
the transitions between the bands and we have observed it in our experiments 
consistently for different bandwidths.

\begin{figure}[t]
  \centering
  \subfloat[250 eigs (Vase-Lion, 200k $V$) ]{\includegraphics[width=0.495\linewidth]{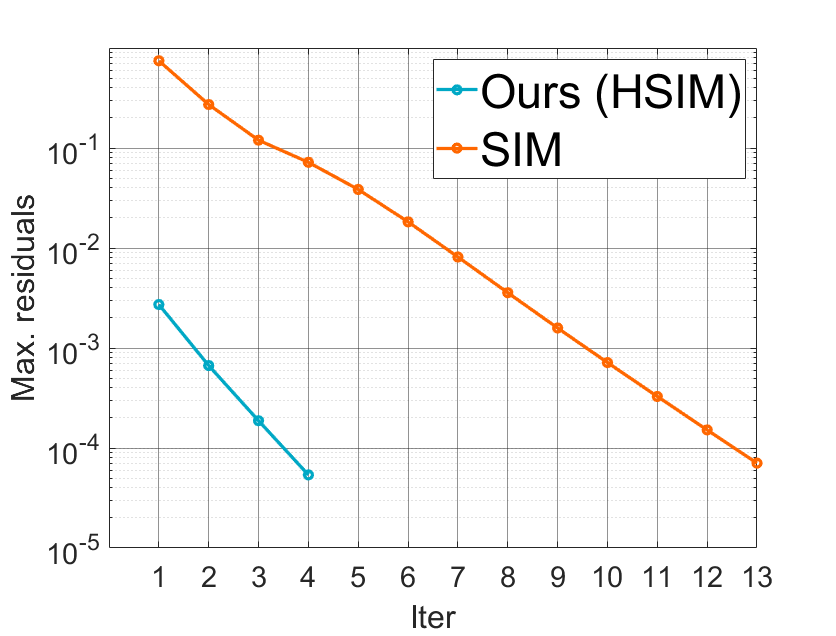}\label{fig:convergence2-1}}
  \hfill
  \subfloat[100 eigs (Eros, 475k $V$)]{\includegraphics[width=0.495\linewidth]{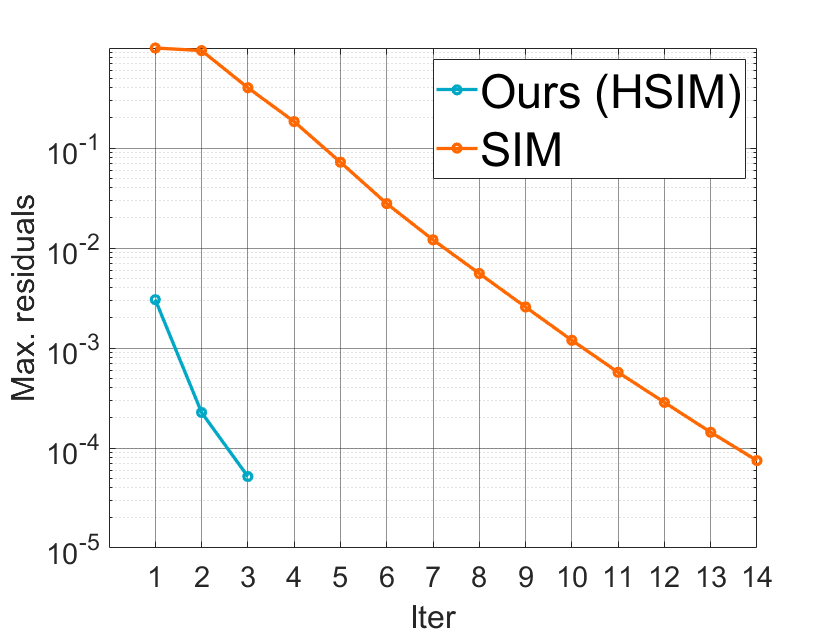}\label{fig:convergence2-2}}  
  \parbox[t]{1.0\columnwidth}{\relax  }
  \caption{\label{fig:residue_comparison}
	Plot of the maximum residual and the numbers of iterations for SIM and HSIM are shown. For HSIM the number of iterations on the finest level is used.}
\end{figure}

\paragraph{SIM}

In addition to the runtimes for Lanczos schemes, Table~\ref{tab:vsAll} also lists times and iteration counts for the (non-hierarchical) SIM. If one compares the number of iterations required by HSIM on the finest level with the number of iterations required by SIM, one sees that HSIM  effectively reduces the number of iterations from 7-8 to 1. Accordingly, we observe that HSIM is 4-8 times faster than SIM. The table lists additional runtimes for an optimized SIM implementation in which the linear systems in step 4 of Algorithm~\ref{alg:basic_SIM} are solved in parallel using OpenMP and the dense eigenproblems, step 7 of Algorithm~\ref{alg:basic_SIM}, are solved on the GPU using \textsc{CUDA}'s \textsc{cuSolver} library.
Figure~\ref{fig:residue_comparison} shows for two examples how the number of iterations changes if a lower convergence tolerance is requested. It can be seen that the increase of iterations is lower for HSIM than for SIM.

\revised{
\paragraph{LOBPCG}
The last column of Table~\ref{tab:vsAll} lists timings for the Locally Optimal Block Preconditioned Conjugate 
Gradient Method (LOBPCG). To generate the timings, we used the author's implementation~\cite{Knyazev2007}. 
\minorrev{
We experimented with Jacobi preconditioners, incomplete Cholesky factorizations and the preconditioner $S-\nu {Id}$, which is suggested in \cite{knyazev2001toward}. The latter produced the best results, which we report. 
} Here $S$ is the stiffness matrix (of the transformed ordinary eigenvalue problem that we also used for the Lanczos solvers), $\nu\in\mathbb{R}$ is approximately in the middle of the first ten eigenvalues \cite{knyazev2001toward}, and $Id$ is the identity matrix.
}
\minorrev{
The results demonstrate that HSIM can solve the eigenproblems faster than LOBPCG with the preconditioners we tested.
}

\begin{figure}[t]
  \centering
  \centering
  \subfloat[32 eigenpairs]{\includegraphics[width=0.49\linewidth]{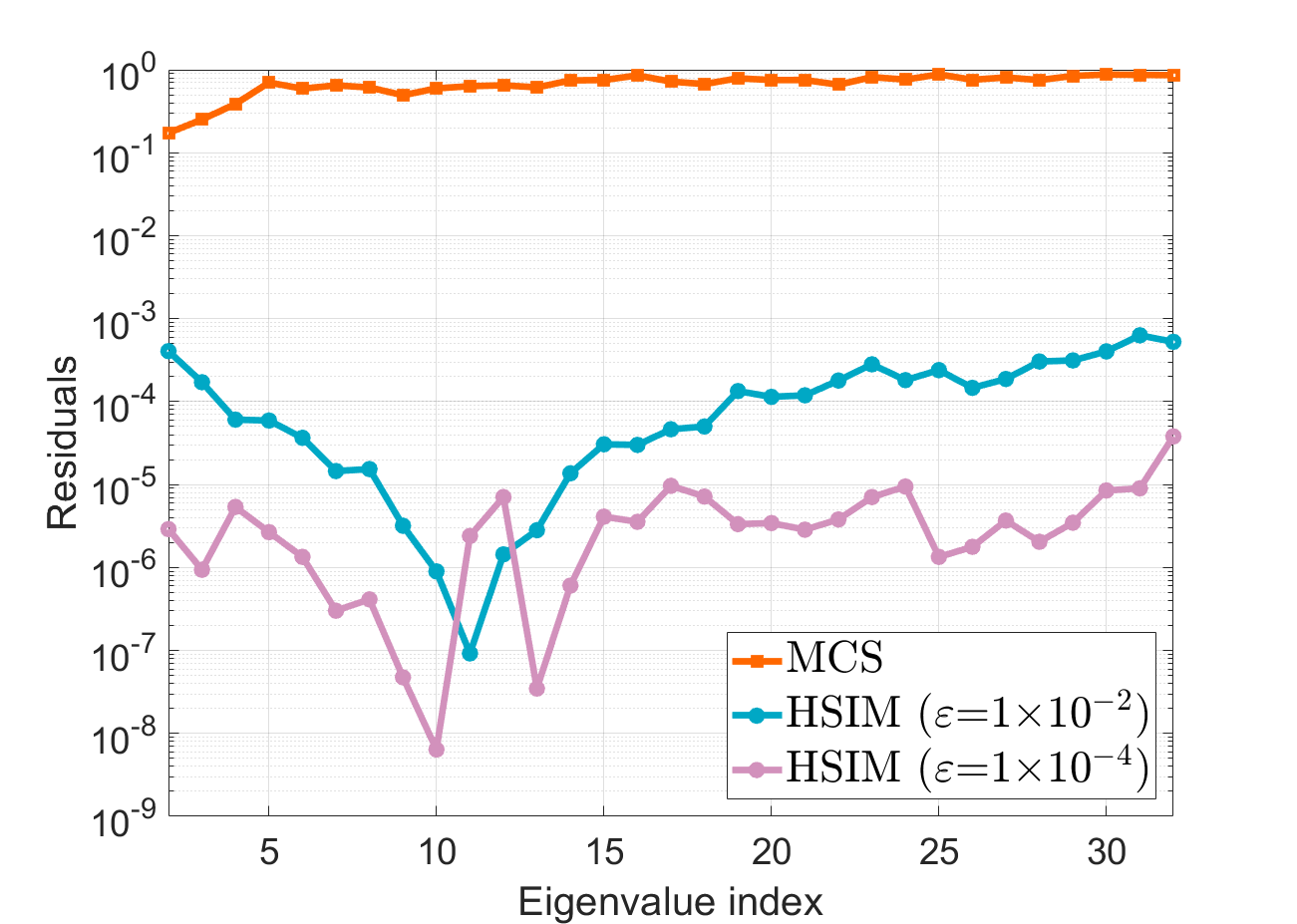}\label{fig:vsBMG1}}
  \hfill
  \subfloat[250 eigenpairs]{\includegraphics[width=0.49\linewidth]{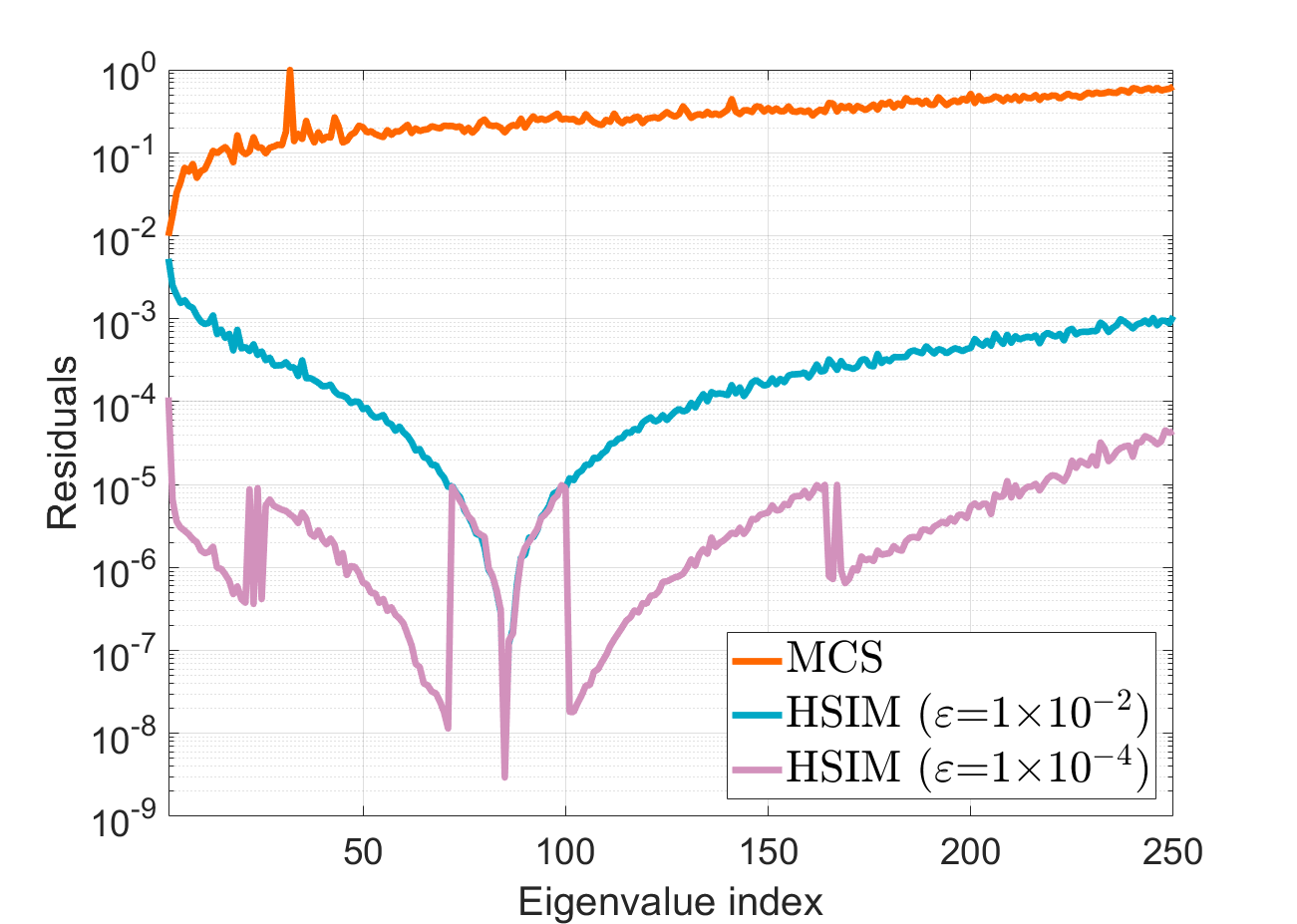}\label{fig:vsBMG2}}  
  \parbox[t]{1.0\columnwidth}{\relax  }
  \caption{
  \label{fig:vsBMG}
   Plot of the residuals of MCS and our novel HSIM eigensolver. 
   Not only HSIM is substantially more accurate, it also has explicit control on the accuracy of the eigenvalues and corresponding eigenvectors.
   (Dragon model, 150k vertices.)
}
\end{figure}
\paragraph{Multilevel correction scheme}
We compare with the multilevel correction scheme (MCS) from \cite{Lin2015,Chen2016}, which 
is an extension of the two-grid scheme from \cite{Hu2011}.  This method has in 
common with our HSIM method that for initialization, an eigenvalue problem on 
the coarsest grid is solved. However, the multilevel iterations differ 
substantially from HSIM. In their method, the coarse space is used in all levels 
and it is enriched by vectors that are computed in the multilevel iterations. An 
essential difference to HSIM is that HSIM reduces the error on each level to the 
desired tolerance margin, while MCS does not offer direct control 
over the accuracy of the solution. The accuracy depends on the approximation 
quality of the coarse grid and the growth rate between the grids. Therefore, an 
aggressive growth rate, which is essential to the performance of our scheme, 
would lead to an increase in approximation error. Another substantial difference 
is that MCS is focused on computing only one or a few eigenpairs. This contrasts this work from our 
setting in which we compute more than a thousand eigenpairs. In Figure~\ref{fig:vsBMG}, 
we show a plot of the accuracy of the eigenpairs computed with MCS 
and HSIM with convergence tolerance $10^{-2}$ and $10^{-4}$. The error produced by MCS 
is orders of magnitude higher than that produced by HSIM. Moreover, the plot on 
the right shows that the error increases with the index of the eigenvalue. This 
illustrates the point that MCS is focused on the computation of a few of the 
lowest eigenpairs. Since the MCS scheme is formulated for regular grids, we use 
our hierarchy with three levels in the comparisons for both schemes, MCS and HSIM.

\section{Conclusion}
We introduce HSIM, a hierarchical solver for sparse eigenvalue problems and evaluate HSIM on the computation of the lowest $p$ eigenpairs of the discrete Laplace--Beltrami operator on triangle surface meshes. 
HSIM first constructs a hierarchy of nested subspaces of the space functions on the mesh. Then, it iterates from coarse to fine over the hierarchy solving the eigenproblem on all levels to the desired accuracy. HSIM is initialized with the solution of the eigenproblem on the coarsest level, which is computed by solving a low-dimensional dense eigenproblem.
Our comparisons show that HSIM outperforms state-of-the-art Lanczos solvers and demonstrate the advantages of the hierarchical approach over the plain SIM.  

We think that the benefits of HSIM over Lanczos and SIM solvers make HSIM attractive for methods in shape analysis and mesh processing. Therefore, we plan to release our implementation of HSIM.

\minorrev{
\paragraph{Future work}
One direction of future work is to explore alternative hierarchies, e.g. operator-dependent bases or wavelets on surfaces. This could improve the performance of HSIM for certain types of operators, such as strongly anisotropic operators. 
Another direction could be to extend the method such that not only the lowest but arbitrary eigenpairs can be efficiently computed. 
Moreover, the method could be improved by further exploring the possibilities of parallelization of the method and by integrating out-of-core techniques for the computation of large eigenbases.
Another aspect is that for a certain complexity of the meshes, the direct solvers will no longer be the most efficient solvers. Then, hierarchical solvers could to be used for the linear systems that need to solved in every iteration. For such an approach, it could be interesting to combine the hierachies used for HSIM and for the linear solves.

A benefit of HSIM is that it directly works for generalized eigenvalue problems, such as (\ref{eq.eigenProb}), and does not require to transform these to ordinary eigenvalue problems. This could be helpful when using the method for solving eigenproblems in which the mass matrix $M$ is not a diagonal matrix, such as the discretization of the Laplace--Beltrami operator with higher-order elements~\cite{Reuter2006}. 
}

\section*{Acknowledgements}
We thank the anonymous reviewers for their constructive feedback. This project is partly supported by the Indonesia Endowment Fund for Education (LPDP) through a doctoral scholarship for Ahmad Nasikun. 
For our experiments, we used models from the AIM@SHAPE repository, the Stanford Computer Graphics Laboratory (Stanford 3D Scanning Repository), Turbosquid, Al-Badri and Nelles (Nefertiti), and INRIA.

\bibliographystyle{ACM-Reference-Format}
\bibliography{Multilevel-Eigensolver}

\appendix
%


%% file: tables/MLsolver.tex
\begin{table}[t]
\revised{
{
	\begin{tabular}{|l|r||r|r|r|r|}
	\hline
	\multirow{2}{*}{\shortstack[l]{Model \\(\#Verts)}}& \multirow{2}{*}{\#Eigs} 			& \multirow{2}{*}{\#Iters} 	& \multicolumn{3}{c|}{Timings of HSIM}  	\\ \cline{4-6} 
									&                         			&                        	& Hier.         		& Solver      	& Total        	\\ \hline \hline	
									& 50                            & F|1			&	{1.9}  				& {  6.2}      	&    8.1   			\\ \cline{2-6} 
									& 250                           & F|1|1			&	{3.7}  				& { 28.4}      	&   32.1   			\\ \cline{2-6} 
	\multirow{-3}{*}{\shortstack[l]{Kitten \\(137k)}}  & 1000          & F|2|1			&	{4.3}  				& {118.9}      	&  123.2   			\\ \hline
                                    & 50                            & F|1           & 3.2                   & 5.3           & 8.5               \\ \cline{2-6} 
                                    & 250                           & F|2|1         & 7.3                   & 41.2          & 48.5              \\ \cline{2-6} 
   \multirow{-3}{*}{\shortstack[l]{Vase-Lion \\(200k)}} & 1000      & F|3|1         & 9.2                   & 188.0         & 197.2             \\ \hline
                                    & 50                            & F|1           & 9.9                   & 28.4          & 38.3              \\ \cline{2-6} 
                                    & 250                           & F|2|1         & 29.1                  & 131.0         & 160.1             \\ \cline{2-6} 
   \multirow{-3}{*}{\shortstack[l]{Knot-Stars\\(450k)}}& 1000       & F|3|1         & 36.3                  & 505.7         & 542.0             \\ \hline
 									& 50                            & F|1			&	{9.2} 				& { 31.9}      	&   41.1   			\\ \cline{2-6} 
									& 250                           & F|2|1			&	{31.6}  			& {122.6}      	&  154.2   			\\ \cline{2-6} 
	\multirow{-3}{*}{\shortstack[l]{Oilpump \\(570k)}} & 1000          & F|3|1			&	{40.3} 				& {650.6}      	&  690.9   			\\ \hline
	                                & 50                            & F|1           & 10.2                  & 65.5          & 75.7              \\ \cline{2-6} 
                                    & 250                           & F|2|1         & 40.6                  & 199.3         & 239.9             \\ \cline{2-6} 
 \multirow{-3}{*}{\shortstack[l]{Red-Circular\\(700k)}} & 1000      & F|4|1         & 55.0                  & 1061.2         & 1116.2             \\ \hline
	\end{tabular}
	}
}
\caption{
\revised{
	Timings of HSIM for the computation of the lowest eigenpairs of the Laplace--Beltrami operator on surface meshes with different numbers of vertices. The error tolerance $\varepsilon$ is set to $10^{-2}$. Individual timings for constructing the hierarchy and for solving the eigenproblem using the hierarchy are listed (in seconds). Meshes are shown in Figure~\ref{fig:eigenfunction}.
	}
}
\label{tab:construction}
\end{table}

%% file: tables/vsAll.tex
\begin{table*}[t]
{
\revised{	\begin{tabular}{|l|r||r|r|r||r|r||r|r|r||r|}
\hline
\multicolumn{1}{|c|}{}                                 & \multicolumn{1}{c||}{}                         & \multicolumn{2}{c|}{SIM}                                 & \multicolumn{1}{c||}{par. SIM} & \multicolumn{2}{c||}{HSIM}                                    & \multicolumn{3}{c||}{Lanczos methods}                                                 & \multicolumn{1}{c|}{Prec. Solver} \\ \cline{3-11} 
\multicolumn{1}{|c|}{\multirow{-2}{*}{Model (\#Vert)}} & \multicolumn{1}{c||}{\multirow{-2}{*}{\#Eigs}} & \multicolumn{1}{c|}{\#Iters} & \multicolumn{1}{c|}{Time} & \multicolumn{1}{c||}{Time}     & \multicolumn{1}{c|}{\#Iters} & \multicolumn{1}{c||}{Time}     & \multicolumn{1}{c|}{\textsc{Matlab}} & \multicolumn{1}{c|}{MH} & \multicolumn{1}{c||}{\textsc{SpectrA}} & \multicolumn{1}{c|}{LOBPCG}       \\ \hline \hline
                                                       & 50                                            & 7                            & 49.4                      & 23.8                          & {\color[HTML]{3531FF} F|1}   & {\color[HTML]{3531FF} 11.7}   & 16.2                        & 27.2                    & 24.4                         & 48.0                              \\ \cline{2-11} 
                                                       & 250                                           & 7                            & 274.5                     & 155.7                         & {\color[HTML]{3531FF} F|2|1} & {\color[HTML]{3531FF} 51.3}   & 94.3                        & 285.3                   & 124.8                        & 268.3                             \\ \cline{2-11} 
                                                       & 1000                                          & 7                            & 1088.0                    & 642.2                         & {\color[HTML]{3531FF} F|2|1} & {\color[HTML]{3531FF} 165.6}  & 921.8                       & 1132.2                  & 1235.5                       & 2601.1                            \\ \cline{2-11} 
                                                       & 2500                                          & 7                            & 3228.2                    & 1930.7                        & {\color[HTML]{3531FF} F|2|1} & {\color[HTML]{3531FF} 529.8}  & 7784.5                      & 2987.1                  & 7552.7                       &  Mem. bound                                 \\ \cline{2-11} 
\multirow{-5}{*}{\parbox{2cm}{Sphere\\(160k)}}                        & 4000                                          & 8                            & 10687.8                   & 8913.0                        & {\color[HTML]{3531FF} F|2|1} & {\color[HTML]{3531FF} 1431.2} & 11745.1                     & 5836.1                  & 13100.1                      &  Mem. bound                                 \\ \hline \hline
                                                       & 50                                            & 8                            & 74.9                      & 42.2                          & {\color[HTML]{3531FF} F|1}   & {\color[HTML]{3531FF} 14.6}   & 18.6                        & 26.4                    & 25.8                         & 126.5                              \\ \cline{2-11} 
                                                       & 250                                           & 8                            & 541.7                     & 300.1                         & {\color[HTML]{3531FF} F|2|1} & {\color[HTML]{3531FF} 79.6}   & 130.3                       & 178.7                   & 185.3                        & 711.5                             \\ \cline{2-11} 
                                                       & 1000                                          & 8                            & 2118.9                    & 1228.0                        & {\color[HTML]{3531FF} F|2|1} & {\color[HTML]{3531FF} 342.1}  & 1549.0                      & 696.4                   & 1359.4                       & 4014.5                            \\ \cline{2-11} 
\multirow{-4}{*}{\parbox{2cm}{Rocker Arm\\(270k)}}                    & 2500                                          & 8                            & 10278.5                   & 8658.4                        & {\color[HTML]{3531FF} F|2|1} & {\color[HTML]{3531FF} 1108.1} & 13018.3                     & 1798.9                  & 9543.0                       & Mem. bound                                   \\ \hline \hline
                                                       & 50                                            & 7                            & 212.5                     & 102.5                         & {\color[HTML]{3531FF} F|1}   & {\color[HTML]{3531FF} 57.9}   & 62.7                        & 100.1                   & 77.7                         & 384.2                             \\ \cline{2-11} 
                                                       & 250                                           & 7                            & 1308.8                    & 664.4                         & {\color[HTML]{3531FF} F|2|1} & {\color[HTML]{3531FF} 206.6}  & 362.5                       & 773.3                   & 675.4                        & 1885.2                            \\ \cline{2-11} 
\multirow{-3}{*}{\parbox{2cm}{Rolling stage\\(660k)}}                 & 1000                                          & 7                            & 8058.2                    & 5358.9                        & {\color[HTML]{3531FF} F|3|1} & {\color[HTML]{3531FF} 937.8}  & 4072.6                      & 3034.5                  & 8396.0                       & Mem. bound                                  \\ \hline
\end{tabular}
}
}
\caption{\revised{Comparison of HSIM to  the (non-hierarchical) SIM, different Lanczos solvers, and LOBPCG. Runtimes are listed in seconds.}}
	\label{tab:vsAll}
\end{table*}

%% file: tables/levels.tex
\begin{table}[t]
\resizebox{1\linewidth}{!} 
{
	\begin{tabular}{|l|r|r|l|r|l|r|l|r|l|r|}
	\hline
	\multicolumn{1}{|c|}{\multirow{2}{*}{\parbox{1.2cm}{Model (\#Verts)}}} 		& \multicolumn{1}{c|}{\multirow{2}{*}{\#Eigs}} & \multicolumn{1}{c|}{\multirow{2}{*}{Tol}}    & \multicolumn{2}{c|}{Level=2}                            & \multicolumn{2}{c|}{Level=3}                            & \multicolumn{2}{c|}{Level=4}                            & \multicolumn{2}{c|}{Level=5}                                                              \\ \cline{4-11} 
	\multicolumn{1}{|c|}{}                       				& \multicolumn{1}{c|}{}                      & \multicolumn{1}{c|}{}                        & \multicolumn{1}{c|}{\#Iters} 	& \multicolumn{1}{c|}{Time}  & \multicolumn{1}{c|}{\#Iters} & \multicolumn{1}{c|}{Time}  & \multicolumn{1}{c|}{\#Iters} & \multicolumn{1}{c|}{Time}  & \multicolumn{1}{c|}{\#Iters} & \multicolumn{1}{c|}{Time} \\ \hline
	\multirow{4}{*}{\parbox{1.2cm}{Rocker Arm (270k)}}          & \multirow{2}{*}{50}		                 & {1e-2}               						& \color[HTML]{0000cd}F|1       	& \color[HTML]{0000cd}{17.2}            & F|1|1                      			& 31.1                       		& F|1|1|1                    & 66.0                       & F|1|1|1|1                  & 123.8                      \\ \cline{3-11} 
																&                                            & {1e-4}               						& \color[HTML]{0000cd}F|3       	& \color[HTML]{0000cd}26.0              & F|3|2   								& {38.2}    						& F|2|2|3                    & 76.8                       & F|2|2|2|2                  & 135.3                      \\ \cline{2-11} 
																& \multirow{2}{*}{2000}                      & {1e-2}               						& F|3       						& 1078.6              					& \color[HTML]{0000cd}{F|2|1}   		& \color[HTML]{0000cd}{650.7}   	& F|1|1|1                    & 885.3                      & F|1|1|1|1                  & 1412.0                     \\ \cline{3-11} 
																&                                            & {1e-4}               						& F|8                        		& 2595.4                      			& \color[HTML]{0000cd}{F|7|4}   		& \color[HTML]{0000cd}{2137.1}  	& F|6|4|4                    & 2967.8                     & F|5|3|3|4                  & 3586.4                     \\ \hline																
	\multirow{3}{*}{\parbox{1.2cm}{Ramses (820k)}}              & 50                                         & \multirow{3}{*}{1e-2}            			& \color[HTML]{0000cd}{F|1}     	& \color[HTML]{0000cd}{45.3}            & F|1|1                      			& 101.0                       		& F|1|1|1                    & 244.0                      & F|1|1|1|1                  & 495.0                      \\ \cline{2-2} \cline{4-11} 
																& 300                                        &                                              & F|2                        		& 246.4                      			& \color[HTML]{0000cd}{F|2|1}   		& \color[HTML]{0000cd}{234.5}   	& F|2|1|1                    & 417.6                      & F|1|1|1|1                  & 700.6                      \\ \cline{2-2} \cline{4-11} 
																& 750                                        &                                              & F|2                        		& 969.6                      			& \color[HTML]{0000cd}{F|2|1}   		& \color[HTML]{0000cd}{657.5}   	& F|2|1|2                    & 1521.1                     & F|2|1|1|1                  & 1607.0                     \\ \hline
	\end{tabular}
}
	\caption{Performance of HSIM with different numbers of levels.
	}
	\label{tab:level}
\end{table}

%% file: tables/support.tex
\begin{table}[b]
\resizebox{\linewidth}{!} 
{
	\begin{tabular}{|l|r||r|r|r|r|r|r|r|r|r|r|}
	\hline
	\multicolumn{1}{|c|}{}                        			& \multicolumn{1}{c||}{}                         & \multicolumn{2}{c|}{2.5}                                       & \multicolumn{2}{c|}{5}                                     & \multicolumn{2}{c|}{7}                                      & \multicolumn{2}{c|}{10}                                    & \multicolumn{2}{c|}{20}                                      \\ \cline{3-12} 
	\multicolumn{1}{|p{0.5cm}|}{\multirow{-2}{*}{\parbox{1.2cm}{Model (\#Verts)}}} 	& \multicolumn{1}{p{0.6cm}||}{\multirow{-2}{*}{\#Eigs}} & \multicolumn{1}{p{0.4cm}|}{\#Iters} & \multicolumn{1}{p{0.13cm}|}{Time}    					& \multicolumn{1}{p{0.5cm}|}{\#Iters} & \multicolumn{1}{c|}{Time}     & \multicolumn{1}{p{0.5cm}|}{\#Iters} & \multicolumn{1}{c|}{Time}     & \multicolumn{1}{p{0.5cm}|}{\#Iters} & \multicolumn{1}{c|}{Time}     & \multicolumn{1}{p{0.5cm}|}{\#Iters} & \multicolumn{1}{c|}{Time}     \\ \hline \hline
															& 100                                           & {F|2}     & {  15.9} 		& {F|2}    							& {  17.0} 							& \color[HTML]{0000cd}{F|1}       	  & \color[HTML]{0000cd}{ 12.2} 		& {F|1}       & { 13.4} 	& {F|2}       & { 23.6}  \\ \cline{2-12} 
															& 400                                           & {F|3|2}   & {  78.3} 		& \color[HTML]{0000cd}{F|3|1}  		& \color[HTML]{0000cd}{  60.4} 		& {F|3|1}     					      & { 69.6} 							& {F|2|1}     & { 74.3} 	& {F|2|1}     & {110.4}  \\ \cline{2-12} 
	\multirow{-3}{*}{\parbox{1.2cm}{Vase-Lion (200k)}}             & 750                                           & {F|4|2}   & { 175.4} 		& {F|3|2}  							& { 174.8} 							& \color[HTML]{0000cd}{F|3|1}    	  & \color[HTML]{0000cd}{138.9} 		& {F|3|2}     & {219.5} 	& {F|3|2}     & {301.5}  \\ \hline
															& 100                                           & {F|2}     & {  55.7} 		& \color[HTML]{0000cd}{F|1}    		& \color[HTML]{0000cd}{  42.2} 		& {F|1}       						  & { 43.0} 							& {F|2}       & { 65.3} 	& {F|2}       & { 71.7}  \\ \cline{2-12} 
															& 400                                           & {F|2|3}   & { 273.6} 		& {F|2|2}  							& { 231.9} 							& \color[HTML]{0000cd}{F|2|1}     	  & \color[HTML]{0000cd}{169.8} 		& {F|4|1}     & {207.2} 	& {F|4|1}     & {292.7}  \\ \cline{2-12} 
	\multirow{-3}{*}{\parbox{1.2cm}{Eros (475k)}}                  & 750                                           & {F|3|4}   & { 710.2} 		& {F|2|3}  							& { 580.5} 							& \color[HTML]{0000cd}{F|4|1}     	  & \color[HTML]{0000cd}{465.5} 		& {F|4|2}     & {567.8} 	& {F|5|2}     & {758.7}  \\ \hline
	\end{tabular}	
}
\caption{Runtimes and iteration counts for different values of the parameter~$\sigma$ that determines the supports of the functions.}
	\label{tab:support}
\end{table}

%% file: tables/shift.tex
\begin{table}[t]
\resizebox{0.975\linewidth}{!} 
{
	\begin{tabular}{|r|c|r|r|r|r|r|r|r|}
\hline
\multicolumn{1}{|c|}{\multirow{2}{*}{\#Eigs}} & \multirow{2}{*}{Residue} & \multicolumn{7}{c|}{Shift ratio}                                                                                                                                                                  \\ \cline{3-9} 
\multicolumn{1}{|c|}{}                        &                          & \multicolumn{1}{c|}{No shift} & \multicolumn{1}{c|}{0.1} & \multicolumn{1}{c|}{0.2} & \multicolumn{1}{c|}{0.25} & \multicolumn{1}{c|}{1/3} & \multicolumn{1}{c|}{0.4} & \multicolumn{1}{c|}{0.45} \\ \hline
50                                            & \multirow{3}{*}{1e-2}    & F|1                           & F|1                      & F|1                      & F|1                       & \textbf{F|1}             & F|2                      & F|2                       \\ \cline{1-1} \cline{3-9} 
250                                           &                          & F|2|1                         & F|2|1                    & F|2|1                    & F|2|1                     & \textbf{F|2|1}           & F|2|1                    & F|2|1                     \\ \cline{1-1} \cline{3-9} 
1000                                          &                          & F|3|1                         & F|2|1                    & F|2|1                    & F|2|1                     & \textbf{F|2|1}           & F|3|2                    & F|4|2                     \\ \hline
50                                            & \multirow{3}{*}{1e-4}    & F|5                           & F|5                      & F|4                      & F|4                       & \textbf{F|4}             & F|4                      & F|5                       \\ \cline{1-1} \cline{3-9} 
250                                           &                          & F|6|4                         & F|6|4                    & F|5|4                    & F|5|3                     & F|5|3                    & \textbf{F|4|3}           & F|5|4                     \\ \cline{1-1} \cline{3-9} 
1000                                          &                          & F|8|4                         & F|8|4                    & F|7|4                    & F|7|4                     & \textbf{F|6|3}           & F|6|4                    & F|7|5                     \\ \hline
\end{tabular}
	}

	\caption{Iteration counts for different choices of shifting values are shown. Computations are done using the Gargoyle model with 85k vertices.
	}
	\label{tab:shift}
\end{table}

%% file: tables/boundary.tex
\begin{table}[b]
{
\revised{	\begin{tabular}{|l|r|r|r|r|r|}
\hline
\multicolumn{1}{|c|}{\multirow{2}{*}{Boundary}} & \multicolumn{1}{c|}{\multirow{2}{*}{\#Eigs}} & \multicolumn{1}{l|}{\multirow{2}{*}{\#Iters}} & \multicolumn{3}{c|}{Timing}                                                          \\ \cline{4-6} 
\multicolumn{1}{|c|}{}                          & \multicolumn{1}{c|}{}                        & \multicolumn{1}{l|}{}                         & \multicolumn{1}{c|}{Hier.} & \multicolumn{1}{c|}{Solve} & \multicolumn{1}{c|}{Total} \\ \hline \hline
\multirow{2}{*}{Dirichlet}                      & 50                                           & F|2                                           & 5.8                        & 30.0                       & 35.8                       \\ \cline{2-6} 
                                                & 250                                          & F|2|1                                         & 19.3                       & 93.7                       & 113.0                      \\ \hline
\multirow{2}{*}{Neumann}                        & 50                                           & F|2                                           & 5.7                        & 29.8                       & 35.4                       \\ \cline{2-6} 
                                                & 250                                          & F|2|1                                         & 18.9                       & 94.9                       & 113.3                      \\ \hline
\end{tabular}
}}
	\caption{\revised{Timings and iteration counts for solving eigenproblems with boundary conditions on the Julius Caesar model with 370k vertices.
	}}
	\label{tab:boundary}
\end{table}